\input amstex
\documentstyle{amsppt}
\magnification 1200
\TagsOnRight
\NoBlackBoxes
\hsize=6.80 truein
\hcorrection{.001in}
\vsize=8.5 truein
\def\Box{\sqcup\!\!\!\!\sqcap}

\define\RR{{R\!\!\!\!\!I~}}

\define\RN{\RR^N}
\define\R+{\RR_+}
\define\del{{\partial}}
\define\eps{\epsilon}
\define\UU{\Cal U}
\define\C{\Cal C}
\define\EE{\Cal E}
\define\sgn{\hbox{\rm sgn}}
\define\supp{\hbox{\rm supp }}
\define\Prob{\hbox{\rm Prob}}
\define\id{\hbox{\rm id}}
\define\diag{\hbox{\rm diag}}
\define \refAmadori            {1} 
\define \refAmadoriColombo     {2} 
\define \refBall               {3}
\define \refBardosLerouxNedelec{4} 
\define \refBenabdallah        {5}
\define \refBenabdallahSerreone{6}
\define \refBenabdallahSerretwo{7}
\define \refBourdelDelormeMazet{8} 
\define \refChenLeFloch        {9}
\define \refCourantFriedrichs  {10}
\define \refDafermos           {11} 
\define \refDiPernaone         {12}
\define \refDiPernatwo         {13}
\define \refDuboisLeFlochone   {14} 
\define \refDuboisLeFlochtwo   {15} 
\define \refDuboisLeFlochthree {16} 
\define \refDubrocaGallice     {17}
\define \refGisclonone   {18} 
\define \refGisclontwo   {19} 
\define \refGisclonSerre {20} 
\define \refGlimm        {21}
\define \refGoodman      {22}
\define \refHartman      {23}
\define \refJosephone    {24} 
\define \refJosephtwo    {25}
\define \refJosephLeFloch{26}
\define \refJosephVeerappa{27} 
\define \refKreiss       {28}
\define \refLaxone       {29}
\define \refLaxtwo       {30}
\define \refLaxthree     {31}
\define \refLeFloch      {32}
\define \refLeFlochtwo   {33}
\define \refLeFlochNedelec{34}
\define \refLeroux       {35} 
\define \refLiYu         {36} 
\define \refLiuone       {37}
\define \refLiutwo       {38}
\define \refLiuthree     {39}
\define \refLiufour      {40}
\define \refMajdaPego    {41} 
\define \refSableTougeronone{42}
\define \refSableTougerontwo{43}
\define \refSerre        {44} 
\define \refSmoller      {45}
\define \refSchwartz     {46} 
\define \refSzepessy     {47} 
\define \refTemple       {48} 
\define \refXin          {49}
\topmatter
This paper was published as : 
K.T. Joseph and P.G. LeFloch, Boundary layers in weak solutions to hyperbolic conservation laws, 
{\it Arch. Rational Mech Anal.} {\bf 147} (1999), 47--88.
 
 \
 
 \
 
\title
BOUNDARY~LAYERS~IN~WEAK~SOLUTIONS\\
OF\\
HYPERBOLIC CONSERVATION LAWS
\endtitle 
\author
K.T. Joseph
\footnote"$^1$"{ \, School of Mathematics,
Tata Institute of Fundamental Research, Homi Bhabha Road,
Bombay 400005, India. E-mail: ktj\@math.tifr.res.in. \hskip9.cm\hfill}
and P.G. LeFloch
\footnote"$^2$"{\, FORMER ADDRESS : 
Centre de Math\'ematiques Appliqu\'ees and 
Centre National de la Recherche Scientifique, 
CNRS UA756, Ecole Polytechnique, 91128 Palaiseau, France.  
NEW ADDRESS (from 2004): 
Laboratoire Jacques-Louis Lions \& Centre National de la Recherche Scientifique, 
Universit\'e de Paris 6, Place Jussieu, 75252 Paris, France. E-mail: {\it LeFloch\@ann.jussieu.fr}
\newline\indent
{\it Key words and phrases\/}:
conservation law, shock wave, boundary layer, vanishing viscosity method, 
difference scheme. \hskip2.cm\hfill}
\endauthor
\rightheadtext{K.T. Joseph and P.G. LeFloch}
\leftheadtext{Boundary layers in weak solutions}
\abstract
This paper is concerned with the initial-boundary value problem
for a nonlinear hyperbolic system of conservation laws. 
We study the boundary layers that may arise in 
approximations of entropy discontinuous solutions. 
We consider both the vanishing viscosity method and 
finite difference schemes (Lax-Friedrichs type schemes, Godunov scheme). 
We demonstrate that {\it different\/} regularization methods 
generate {\it different\/} boundary layers. 
Hence, the boundary condition can be formulated only if 
an approximation scheme is selected first.   
Assuming solely uniform $L^\infty$ bounds on the approximate 
solutions and so dealing with $L^\infty$ solutions, 
we derive several entropy inequalities satisfied by the boundary layer
in each case under consideration. A Young measure is   
introduced to describe the boundary trace. When a uniform bound on 
the total variation is available, the boundary Young measure reduces to 
a Dirac mass. Form the above analysis, we deduce 
several formulations for the boundary condition 
which apply whether the boundary is characteristic or not. 
Each formulation is based a {\it set of admissible boundary
values\/}, following Dubois and LeFloch's terminology in 
``Boundary conditions for nonlinear hyperbolic 
systems of conservation laws'', J. Diff. Equa. 71 (1988), 93--122. 
The local structure of those sets and the well-posedness of the 
corresponding initial-boundary value problem are investigated. 
The results are illustrated with convex and nonconvex conservation 
laws and examples from continuum mechanics. 
\endabstract
\endtopmatter
\document

\subheading {1. Introduction}
 
This paper considers the initial-boundary value problem for an 
hyperbolic system of conservation laws
$$
\del_t u + \del_x f(u) = 0, \qquad u(x,t)\in \UU \subset \RN, 
\quad x>0, \, t > 0, 
\tag 1.1
$$
supplemented with 
\roster
\item an initial condition at time $t=0$
$$
u(x,0) = u_I(x), \qquad x>0, 
\tag 1.2
$$
\item the entropy inequality
$$
\del_t U(u) + \del_x F(u) \le 0,
\tag 1.3 
$$
\item and a weak form of the following Dirichlet
boundary condition at $x=0$
$$
u(0,t) = u_B(t), \qquad t>0.
\tag 1.4 
$$ 
\endroster

Indeed the hyperbolic problem (1.1)--(1.4) 
is usually not well-posed when the boundary data is required to be
assumed in the (strong) sense (1.4), even when (1.1) 
is a linear system (cf.~Kreiss \cite{\refKreiss}). 
It is the objective of this paper to provide a general framework
which leads to (mathematically correct) formulations for the boundary condition. 
Following Dubois-LeFloch 
\cite{\refDuboisLeFlochtwo}, our strategy
is to reformulate (1.4) in the (weak) form
$$
u(0+,t) \in \EE(u_B(t)), \qquad t>0,
\tag 1.5
$$
where $\EE(u_B(t)) \subset \UU$ is a time-dependent set
(the set of {\it admissible boundary values\/})
to be defined from the boundary data, 
and $u(0+,t)$ is the trace (its existence is discussed in this paper)
of the solution $u$ at the boundary.  
We shall consider several methods of 
approximation for the problem (1.1)--(1.4), 
including the artificial vanishing viscosity method 
and a class of finite difference schemes, for which 
the boundary condition (1.4) can be easily implemented.  
As the approximation parameter goes to zero, 
a sharp transition layer generally develops near the boundary $\big\{x=0\big\}$  
and the limiting solution does not satisfy the boundary condition (1.4). 
Our aim in this paper is to provide some contribution to the following 
program: 
perform a rigorous analysis of the boundary layer for weak solutions, 
then derive several suitable definitions for the set in (1.5), and finally 
investigate the structure of the latter to decide whether the 
boundary-value problem is well-posed.

In (1.1), $\UU$ is assumed to be a convex and open subset of $\RN$, 
the flux-function $f: \UU \to \RN$ to be a smooth mapping,  
and the initial data $u_I$ to belong to 
$L^\infty(\RR_+, \UU)$. It will be convenient to assume that
the boundary data $u_B$ has bounded total variation on any interval $[0,T]$ 
for all $T>0$.
It is assumed that (1.1) admits at least one strictly convex
entropy pair. By definition, a pair of functions 
$(U,F): \UU \to \RR\times\RR$ of class $\C^2$ is called a convex (or strictly
convex) entropy pair iff $\nabla F^T = \nabla U^T \nabla f$
and the Hessian matrix $\nabla^2 U$ is non-negative (or positive definite).
The existence of at least one strictly convex
entropy pair implies that (1.1) is hyperbolic. 
For background on hyperbolic systems, we refer 
to Lax \cite{\refLaxone, \refLaxtwo, \refLaxthree}, 
Dafermos \cite{\refDafermos} and Smoller \cite{\refSmoller}, concerning
the theory of existence of entropy solutions to the pure Cauchy problem,
to Glimm \cite{\refGlimm} and Liu \cite{\refLiufour} 
for initial data with small total variation, 
and DiPerna \cite{\refDiPernaone,\refDiPernatwo} 
for systems of two equations with $L^\infty$ initial data.

This paper contributes to establishing a framework for the 
initial-boundary value problem for (1.1). 
It is intended to pursue the efforts initiated in recent years on this problem 
(Cf.~review below). In particular we built upon the recent contributions in 
Gisclon-Serre \cite{\refGisclonSerre} and Xin \cite{\refXin}, who
studied the boundary layers associated with the vanishing viscosity 
approximations assuming 
the solution to the hyperbolic problem be smooth. 
A formal asymptotic expansion is introduced in
\cite{\refGisclonSerre, \refXin} and the convergence including 
$L^2$ error estimates is proven for the boundary layer
in the smooth regime.
For linear hyperbolic systems, Joseph \cite{\refJosephtwo}
constructed boundary layers explicitly and obtained error estimates in
$L^2$ Sobolev space.

One of the motivations here is to treat several approximation methods
simultaneously and compare the results obtained with each of them. 
We consider the vanishing viscosity method, a class of Lax-Friedrichs 
type schemes, and the Godunov scheme. 

In Section 2, we rigorously derive conditions satisfied by the boundary layer, 
which take the form of a family of {\it boundary entropy inequalities\/}
and a {\it boundary layer equation\/}. 
The regularity of the relevant traces at the boundary are discussed. 
The whole analysis is performed by assuming only a uniform $L^\infty$ 
bound on the approximate solutions; in particular no assumption is required on the 
regularity of the limiting solution to (1.1). 
Since high frequency oscillations in the approximate solutions 
can not be a~priori excluded, 
the conditions above are formulated in terms of a boundary Young measure 
associated with the boundary layer. 
Note that, in the derivation of Section 2, 
the boundary is possibly characteristic, i.e.~the 
eigenvalues of the matrix $\nabla f(u)$ may vanish for certain values
of $u$. 

Observe also that, in general, the equations and inequalities we derive 
depend upon the approximation method in use. 
Fundamentally the boundary condition can not be formulated from the 
mere knowledge of the function $u_B$, but depend upon the underlying 
``physical'' regularization. 
This feature arises in weak solutions to many nonlinear hyperbolic
problems. See, for instance, the review paper by LeFloch
\cite{\refLeFlochtwo} on regularization-sensitive shock waves. 

In Section 3, we introduce several sets of admissible boundary values 
and investigate their local structure. 
When the boundary is non-characteristic, we establish that the sets based 
on the boundary layer equations are manifold with the ``correct'' dimension.  
That is, the corresponding initial-boundary value problem is well-posed,
at least for constant boundary and initial data 
(a generalization to the Riemann problem). We also prove a similar (but
stronger) result for the set
based on the boundary layer equation derived by the Godunov scheme. 
Strictly speaking this scheme does not produce
any boundary layer; however analyzing that scheme
 leads to a formulation of the boundary condition as it was first pointed 
out in \cite{\refDuboisLeFlochtwo, \refDuboisLeFlochthree}. 
We recall that setting the boundary condition via an upwinding
difference scheme is a classical idea in the computing literature.

Sections 4 is ddevoted to studying several examples of particular interest. 
It is expected that, in general, different approximation method 
for (1.1) leads to a different set in (1.5). However we prove in Section 4, 
for both convex and non-convex conservation laws, that this
is not the case when $N=1$. In other words the boundary layer for
the scalar conservation laws is independent of the approximation method. 
The same is true of the linear hyperbolic systems; and 
we conjecture that this also holds 
for the nonlinear systems in the class with coinciding shock and rarefaction curves
introduced by Temple \cite{\refTemple}. We also consider
examples from continuum mechanics, i.e.~the system 
of nonlinear elasticity and the system of gas dynamics. 

To complete this presentation, we give a short 
overview of the literature on the boundary conditions for (1.1). 
Most of the activity was restricted to scalar equations, i.e.~$N=1$. 
The pioneering work by Leroux \cite{\refLeroux} and Bardos-Leroux-Nedelec 
\cite{\refBardosLerouxNedelec}  based on the vanishing viscosity method 
provides a derivation of 
 ``the''correct formulation of the boundary condition for multidimensional
{\it scalar\/} conservation laws.
Specically, \cite{\refBardosLerouxNedelec} shows that 
(1.4) should be replaced by the weaker statement:
$$
\big(\sgn(u(0+,t) - k) - \sgn(u_B(t)-k) \big) \, 
  \bigl(f(u(0+,t)) - f(k)\bigr) \ge 0 
\qquad \text{ for all } k \in \RR, 
\tag 1.6 
$$
where $\sgn (a) = -1$ if $a<0$,  $\sgn (a) = 0$ if $a=0$,
and $\sgn (a) = 1$ if $a>0$. 
The convergence of finite difference schemes,
again for scalar equations, is established by Leroux in an unpublished work:
it is remarkable that the finite difference scheme approach leads
to the same formulation (1.6) of the boundary condition. 
The condition is used by LeFloch \cite{\refLeFloch} 
in order to extend Lax's explicit formula \cite{\refLaxtwo} to the
initial-boundary value problem. Joseph \cite{\refJosephone}
used the vanishing viscosity method and the Hopf-Cole transformation 
to extend Lax's formula for the inviscid Burgers equation.  
Another derivation is given by Joseph and Veerappa Gowda  
\cite{\refJosephVeerappa}; see also Gisclon \cite{\refGisclonone} and
LeFloch-Nedelec \cite{\refLeFlochNedelec}.We also refer to the 
paper \cite{\refSzepessy} by Szepessy for a very general result of 
existence and uniqueness.

The statement (1.6) is a special case  (when applied to Kruzkov entropies) of a more 
general inequality: 
$$
F(u(0,t)) - F(u_B(t)) 
- \nabla U(u_B(t)) \bigl(f(u(0,t)) - f(u_B(t))\bigr) 
\le 0, 
\tag 1.7 
$$
which has to hold for every convex entropy pair $(U,F)$. 
The latter was derived formally using the vanishing viscosity method
in Dubois-LeFloch \cite{\refDuboisLeFlochtwo}, 
who pointed out that (1.7) holds even when $N \ge 2$ 
and introduced the notion of set of admissible boundary values, 
cf.~(1.5). 
These inequalities were obtained independently by Bourdel-Delorme-Mazet 
\cite{\refBourdelDelormeMazet} 
based on an analysis of the characteristics of the system (1.1), 
and by Benabdallah \cite{\refBenabdallah} for a specific system.
The first result of existence
for the initial-boundary value problem
for a system was given by Benabdallah-Serre 
\cite{\refBenabdallahSerreone, \refBenabdallahSerretwo}: 
the vanishing viscosity method
applied to the $p$-system of gas dynamics converges to a 
solution to (1.1) satisfying the set of inequalities (1.7).  

The Glimm scheme with various type of boundary conditions was 
studied by Liu, for instance \cite{\refLiuone, \refLiutwo, \refLiuthree}. 
In the case that the boundary is assumed to be non-characteristic
and the number of boundary conditions is equal to the number of
positive eigenvalues of the matrix $\nabla f$,
Goodman proves the convergence of the Glimm scheme
in his unpublished thesis \cite{\refGoodman}; 
cf.~also Dubroca-Gallice \cite{\refDubrocaGallice} and 
Sabl\'e-Tougeron \cite{\refSableTougeronone, \refSableTougerontwo}

More recently Amadori \cite{\refAmadori, \refAmadoriColombo} 
used the formulation in \cite{\refDuboisLeFlochtwo}
and proved the convergence 
of a front tracking scheme in the characteristic case. 
In particular, Amadori establishes that a condition of the form 
(1.5) can be satisfied pointwise except at countably many times. 

Finally we refer to the IMA report \cite{\refJosephLeFloch} by the authors 
for an extended version of the present article. 

\subheading{2. Boundary Layers in Weak Solutions}

In this section, we consider sequences of approximate solutions to the 
initial boundary value problem (1.1)--(1.4), and aim at
characterizing their limiting behavior near the boundary. 
Here we rigorously 
derive entropy inequalities satisfied by the boundary layer. 
We deal with a sequence of $L^\infty$ functions  
with uniformly bounded amplitude. 
As is well-known, for general systems of conservation laws, proving 
the strong convergence of a sequence of approximate solutions 
is an open problem. 
It seems therefore natural to formulate those 
entropy inequalities in terms of a Young measure (for instance 
Ball \cite{\refBall} for this concept) associated with the 
sequence of approximate solutions.  Further analysis can be
performed on a case by case basis only. 

In the following, certain averages will be shown to belong to the space 
$BV(\RR_+)$ of functions of locally bounded total variation, 
i.e.~measurable and bounded 
functions $w: \RR_+ \to \RR$ whose distributional derivative is a bounded
Borel measure on every interval $(0,T)$ for all $T>0$. 
We denote by $TV_0^T(w)$ the total variation, and by 
$\|w\|_{BV(0,T)} = \|w\|_{L^\infty(0,T)} + TV_0^T(w)$ the norm, of a BV function $w$
on an interval $(0,T)$. By convention, 
a BV function will be always normalized by 
selecting its right continuous representative. 
\vskip.3cm 

\noindent{\bf 2.1 Vanishing Viscosity Method.} 
Let $u^\eps$ be the approximate solutions obtained by solving the following 
parabolic regularization of (1.1)-(1.4): 
$$
\del_t u^\eps + \del_x f(u^\eps) = \eps \, \del^2_{xx} u^\eps, 
\qquad \quad x>0, \, t > 0, 
\tag 2.1
$$
$$
u^\eps(x,0) = u_I^\eps(x), \qquad x>0, 
\tag 2.2
$$
$$
u^\eps(0,t) = u_B^\eps(t), \qquad t>0. 
\tag 2.3 
$$
The smooth functions $u_I^\eps \in L^\infty(\RR_+)$ and $u_B^\eps \in BV(\RR_+)$
are chosen to be uniformly bounded and a.e. convergent
approximations of the corresponding data $u_I$ and $u_B$.  
We assume the existence of a (smooth enough) solution $u^\eps$ to the problem
(2.1)--(2.3). Note that compatibility conditions at $(x,t)=(0,0)$, 
such as $u_I^\eps(0) = u_B^\eps(0)$, are implicitly required. 
We shall also assume that 
$$
u^\eps \text{ is uniformly bounded in } L^\infty(\RR\times\RR_+). 
\tag 2.4 
$$

We introduce a new function $v^\eps$ by setting 
$$
v^\eps(y,t) = u^\eps(\eps y, t), 
\tag 2.5 
$$
so that the system of equations (2.1) transforms into
$$
\eps \, \del_t v^\eps + \del_y f(v^\eps) 
     =  \del^2_{yy} v^\eps.
\tag 2.6
$$
It is expected that
the ($\eps \to 0$) limit of the $v^\eps$'s will give us
a good description of
the boundary layer at $x=0$, at least under additional assumptions, 
although a different scaling may more adapted in certain circumstances.

By definition (e.g.~Ball \cite{\refBall}), 
a Young measure associated with a sequence $u^\eps$ satisfying (2.4)
is a weak-star measurable 
mapping $\nu$ from the $(x,t)$ plane to the space $\Prob(\UU)$ of 
all probability measures (i.e.~non-negative measures with mass one)
with the property that 
for every continuous function $g:\UU \to \RR$  
$$
g(u^\eps) \rightharpoonup <\nu,g> \quad \text{ weakly--$\star$ in } L^\infty(\RR^2_+).
\tag 2.7
$$
In view of (2.4), 
the functions $v^\eps$ also are uniformly bounded in $L^\infty(\RR\times\RR_+))$.
We denote by $\mu$ a Young measure associated with the functions $v^\eps$.

\proclaim{Theorem 2.1} 
The following statements hold for all convex entropy pairs
$(U,F)$ associated with the system $(1.1)$, all functions 
$\theta \in BV(\RR_+)$, and any bounded interval $(T_1, T_2)$. \par 
1) \, When $\theta(t) \ge 0$, the distribution  
$$
y \, \mapsto \, \int_{T_1}^{T_2} <\mu_{y,t}, F> \, \theta(t) \, dt 
- \frac {d}{dy}\int_{T_1}^{T_2} <\mu_{y,t},U> \, \theta(t) \, dt
$$ 
is in fact a function of locally bounded variation and 
thus is defined pointwise as a right continuous function. 
There exists a Young measure $\mu_{0,t}$, 
such that the following limit exists and is given by $\mu_{0,t}$: 
$$
\lim_{y \to 0+} \int_{T_1}^{T_2} <\mu_{y,t},U> \, \theta(t) \, dt
= 
\int_{T_1}^{T_2} <\mu_{0,t},U> \, \theta(t) \, dt. 
$$ 
When $\theta(t) \ge 0$, the function 
$$
x \, \mapsto \,  \int_{T_1}^{T_2} <\nu_{x,t},F> \theta(t) \, dt
$$
has locally bounded variation. There exists 
a Young measure $\nu_{0,t}$, the ``trace'' of $\nu_{x,t}$ at $x=0$, 
such that the following limit exists and is given by $\nu_{0,t}$: 
$$
\lim_{x \to 0+} \int_{T_1}^{T_2} <\nu_{x,t},F> \, \theta(t) \, dt
= 
\int_{T_1}^{T_2} <\nu_{0,t},F> \, \theta(t) \, dt. 
$$ 
When $(U,F) = (u_j, f_j)$, $1 \leq j \leq N$, 
all of the results above still hold when the function 
$\theta$ has no specific sign. 
\par
2) \, 
For all $0<y_1<y_2$ and in the sense of distributions for $t \in \RR_+$, 
one has 
$$
\aligned 
F(u_B) + \nabla U(u_B) \bigl(<\nu_{0,t}, f> - f(u_B)\bigr)
& \ge \, <\mu_{y_1,t}, F> - \del_y <\mu_{y_1,t},U>  \\
& \ge \, <\mu_{y_2,t}, F> - \del_y <\mu_{y_2,t},U> \\
& \ge  \, <\nu_{0,t}, F>. 
\endaligned 
\tag 2.8
$$
3) \, Moreover one has   
$$
\mu_{0,t} = \delta_{u_B(t)} \qquad \text { a.e. } t \in \RR_+
\tag 2.9 
$$
and, when $\theta \ge 0$,
$$
\lim_{y \to \infty} \biggl(\int_{T_1}^{T_2} <\mu_{y,t}, F> \, \theta(t) \, dt 
- \frac {d}{dy}\int_{T_1}^{T_2} <\mu_{y,t},U> \, \theta(t) \, dt \biggr)
\ge \, 
\int_{T_1}^{T_2} <\nu_{0,t}, F> \, \theta(t) \, dt. 
\tag 2.10 
$$ 
$\hfill\Box$
\endproclaim

Theorem 2.1 provides a rigorous basis to the formal asymptotic 
expansion approach. We collect here several important remarks, 
including the property that the Young measures $\nu$ and $\mu$ reduce 
to Dirac masses when a uniform total variation bound is available. 

First of all, 
the inequalities (2.8) actually hold in the (stronger) sense:
$$
\aligned 
&  \int_{T_1}^{T_2} \biggl(F(u_B(t)) + \nabla U(u_B(t)) \bigl(<\nu_{0,t}, f> 
    - f(u_B(t))\bigr)\biggr) \, \theta(t) \, dt \\
&  \ge 
\int_{T_1}^{T_2} <\mu_{y_1,t}, F> \, \theta(t) \, dt 
- \frac {d}{dy}\bigl(\int_{T_1}^{T_2} <\mu_{y,t},U> \,
       \theta(t) \, dt\bigr)_{\big|y = y_1}  \\ 
&  \ge 
\int_{T_1}^{T_2} <\mu_{y_2,t}, F> \, \theta(t) \, dt 
- \frac {d}{dy}\bigl(\int_{T_1}^{T_2} <\mu_{y,t},U> \, \theta(t) \, 
      dt\bigr)_{\big|y = y_2}         \\ 
&  \ge 
\int_{T_1}^{T_2} <\nu_{0,t}, F> \, \theta(t) \, dt 
\endaligned 
$$
for all non-negative $\theta \in BV(\RR_+)$ and all $0<y_1<y_2$. 
Observe that this is a stronger statement 
than the convergence in the sense of distributions 
since $\theta$ is a function of bounded total variation,
not necessarily having compact support in $(T_1, T_2)$, 
rather than a smooth function with compact support. 
All the formulas to be derived in this section hold in this sense.
Note also that (2.10) is an immediate consequence of (2.8) by taking $y \to \infty$.

The following inequalities, rigorously derived in Theorem 2.1, 
$$
F(u_B) + \nabla U(u_B) \bigl(<\nu_0,f> - f(u_B)\bigr) 
\, \ge \, <\nu_0,F> 
\tag 2.11
$$
will be referred to as the {\it boundary entropy inequalities\/}. 
They do not refer explicitly to the boundary layer itself 
but only to its limiting values. 

The inequalities (2.8) also contain constraints for the boundary layer. 
In particular, using the trivial entropies $(U,F) = \pm (u_j, f_j(u))$, 
$1 \leq j \leq N$, 
in (2.8) leads us to the equation 
$$
<\mu, f> - \del_y <\mu,\id> = <\nu_{0,t}, f>,  
\tag 2.12
$$
where the right hand side is independent of the variable $y$ and only depends on $t$. 

For scalar equations and when the method of compensated compactness 
due to Murat-Tartar applies
(i.e., mainly, for systems of two conservation laws), it is known that 
$\nu$ is a Dirac mass concentrated at a point $u(x,t)$ which is an 
entropy weak solution. In those two situations,  
it is conceivable that the Young measure $\mu$ also would be a Dirac mass. 

If one assumes that $\mu$ is a Dirac mass, say 
$$
\mu_{y,t} = \delta_{v(y,t)} \quad \text { for almost every } \, (y,t) 
\tag 2.13
$$
with $v \in L^\infty$, then the formulas in Theorem 2.1 
take a much simpler form.  
Namely if (2.12) holds, then (2.12) becomes what will be referred to as 
{\it boundary layer equation\/}:
$$
f(v) - \del_y v = <\nu_0,f>. 
\tag 2.14
$$
This is nothing but the equation that would be obtained 
{\it formally\/} by plugging an asymptotic expansion of the form 
$u_\eps(x,t) = u(x,t) + v(x/\eps, t) + O(\eps)$ in the equations (2.1). 
More generally, if (2.12) holds, the inequalities (2.8) become
$$
\aligned 
F(u_B) + \nabla U(u_B) \bigl(<\nu_0,f> - f(u_B)\bigr) 
&  \ge \, F(v(y_1))- \del_y U(v)_{|y=y_1}  \\ 
&  \ge \, F(v(y_2))- \del_y U(v)_{|y=y_2} \\
&  \ge \, <\nu_0,F>.
\endaligned 
$$

When $\nu_0$ also is a Dirac mass for a.e. $t$, say $\nu_{0,t} = \delta_{u_0(t)}$, 
for instance when $u$ has bounded variation in $x$ and 
so admits a trace at $x=0$ in a classical sense, then 
the {\it boundary layer equation\/} (2.14) becomes 
$$
f(v) - \del_y v = <\nu_0,f>. 
\tag 2.15
$$
and the {\it boundary entropy inequalities\/} (2.11) take the form 
$$
F(u_0) - F(u_B) - \nabla U(u_B) \big(f(u_0) - f(u_B) \big) \le 0,
\tag 2.16
$$
which was derived in Dubois-LeFloch \cite{\refDuboisLeFlochone, \refDuboisLeFlochtwo} 
by assuming a uniform BV bound on the $u^\eps$.

Note finally that the behavior of $\mu_{y,t}$ as $y \to \infty$
is controled by the set of inequalities (2.10), only.
If it is assumed that $v$ has a limit in a classical sense and $\del_y v(y,t)\to 0$ as 
$y \to \infty$, then we can set
$$
v_\infty(t) \equiv \lim_{y \to \infty} v(y,t)
$$
and (2.10) becomes
$$
F(v_\infty)  \ge F(u_0) \quad \text{ for all entropy flux } \, F 
\tag 2.17
$$
(the flux $F$ must be associated with a convex entropy).  
In fact (2.17) need not imply 
$$
v_\infty(t) = u_0(t). 
\tag 2.17'
$$
However (2.17) does imply
$$
f(v_\infty(t)) = f(u_0(t))  
$$
so, in the non-characteristic case i.e.~when $\nabla f$ is invertible,
(2.17) implies (2.17'). In the characteristic case, (2.17') may very well
be violated. This difficulty is related to the choice of the scaling in the definition
of the functions $v^\eps$. 
Cf.~the examples in Sections 4 and 5. 

Uniform bounds on the total variation of $u^\eps$ 
are available for scalar equations, linear systems and systems in the so-called Temple's class 
having coinciding shock and rarefaction curves. In the general case we have~: 

\proclaim{Corollary 2.2 } \rm Assume that the solutions to the 
boundary-value problem $(2.1)$-$(2.3)$ additionally satisfy the 
bound
$$
TV(u^\eps(t)):= \int_0^\infty \bigl|\del_x u^\eps(t)\bigr| \, dx \leq C, 
$$
where $C$ is independent of $\eps$. Then the Young measures 
$\mu$ and $\nu_0$ in Theorem 2.1 reduce to Dirac masses, i.e., 
$$
\mu_{y,t} = \delta_{v(y,t)} \qquad \text{ and } 
\qquad \nu_{0,t} = \delta_{u_0(t)}, 
$$
and the functions $v=v(y,t)$ and $u_0= u_0(t)$ satisfy the conditions 
$(2.14)$ and $(2.16)$ for almost every $y,t$.  
$\hfill\Box$ 
\endproclaim


\proclaim{Proof of Theorem 2.1} \rm We decompose the proof into several steps. 
For the whole of this proof, we denote by $(U,F)$ a given convex entropy pair. \par 
\noindent{\bf Step 1:} \, Preliminaries. \par 
We gather here several properties of $\nu$ and $\mu$ that are readily obtained. 
Let us multiply the equation (2.6) by the gradient of $U$ and obtain
$$
\aligned 
\eps \, \del_t U(v^\eps) + \del_y \bigl(F(v^\eps) - \del_y U(v^\eps)\bigr) 
     & =   -  \nabla^2U(v^\eps) \cdot\bigl(\del_y v^\eps, \del_y v^\eps\bigr)\\
     & \le 0. 
\endaligned 
\tag 2.18
$$
Using the definition of the Young measure $\mu$, 
it is a simple matter to pass to the limit in the inequality (2.18). 
For any $\theta \in BV$ and uniformly in $y \in \RR_+$, we have
$$
\aligned 
\eps \bigl|\int_{T_1}^{T_2} \del_t U(v^\eps) \, \theta \, dt \bigr| 
\le \, &  \bigl|\int_{T_1}^{T_2} \eps\, U(v^\eps) \, \del_t \theta \, dt \bigr| \,
          + \, \eps \, \bigl|\bigl[U(v^\eps) \, \theta 
                                \bigr]_{T_1}^{T_2}\bigr|  \\
\le \, & O(1) \,  \|\theta\| _{BV} \, \|U(v^\eps)\|_{L^\infty} \, \to \, 0, 
\endaligned 
$$
so we obtain 
$$
\del_y \biggl(\int_{T_1}^{T_2} <\mu_{y,t}, F> \, \theta \, dt 
- \frac {d}{dy} \int_{T_1}^{T_2} <\mu_{y,t},U> \, \theta \, dt \biggr) \, \le \, 0, 
\tag 2.19
$$
which provides the second inequality in (2.8). 
Therefore time-averages of the function 
$<\mu_{y,t}, F> - \del_y <\mu_{y,t},U>$ are non-increasing, and so have bounded variation 
on any compact set. The limits as $y \to 0+$ or $y \to +\infty$ exist, although 
at this stage of the proof, we can not exclude that those limits 
could be $\pm \infty$.
We shall see later that actually $<\mu_{y,t}, F> - \del_y <\mu_{y,t},U> \in L^\infty$. 
Moreover the function 
$$
\int_{T_1}^{T_2} <\mu_{y,t}, U> \theta(t) \, dt 
$$
has a trace at $y= 0$, which defines $<\mu_{0,t}, U>$. 
Note also that (2.19) with the choices $(U,F) = \pm (u_j, f_j)$, $1 \leq j \leq N$,
leads us to
$$
<\mu_{y,t}, f> - \del_y <\mu_{y,t},\id> = C_*(t), 
\tag 2.20 
$$
where $C_*(t)$ has to be determined. In fact it will be immediate from the 
results in Step 5 below that  
$$
C_*(t) \, = \, <\nu_{0,t}, \id> \qquad \text{ for a.e. } \, t>0. 
$$

Similarly, following DiPerna \cite{\refDiPernatwo} and 
using the Young measure $\nu_{x,t}$ associated with $u^\eps$,
one can pass to the limit in (2.1) and obtain the entropy inequality:
$$
\del_t <\nu_{x,t}, U> + \del_x <\nu_{x,t},F> \le 0. 
\tag 2.21 
$$
{}From (2.21), we deduce first that, for any smooth function $\theta(t) \ge 0$, 
$$
\frac{d}{dx} \int_{T_1}^{T_2} <\nu_{x,t}, F> \theta(t) \, dt 
\le
\int_{T_1}^{T_2} <\nu_{x,t}, U> \del_t\theta(t) \, dt  
\le O(1) \, \|\theta\|_{BV}.
\tag 2.22 
$$
For $\theta$ fixed, the right hand side of (2.22) is a constant, thus 
its left hand side is a locally bounded Borel measure and the function 
$$
g_\theta (x) \, \equiv \, \int_{T_1}^{T_2} <\nu_{x,t}, F> \theta(t) \, dt
$$
has bounded total variation. 
Therefore the trace $\nu_{0,t}$ introduced in Theorem 2.1 exists, 
at least on entropy fluxes. 
This gives a meaning to the last term in the right hand side of (2.8). 
In fact it is possible to establish the estimate 
$$
TV( g_\theta) \, \le \, O(1) \, \|\theta\|_{BV}
$$
for arbitrary functions $\theta \in BV$. (For such $\theta$, 
(2.22) can be obtained directly from (2.1).) Thus the trace $\nu_{0,t}$ 
exists for $\theta \in BV$ as well.

Observe that the traces $\mu_{0,t}$ and $\nu_{0,t}$ are uniquely
determined on entropies and entropy fluxes, respectively.
They can be easily extended as Young measures defined on the whole set of
continuous functions, in a non-unique way however. 
Namely, to construct $\mu_{0,t}$, take any sequence $y_k \to 0$ and 
consider a Young measure associated with the sequence of measures 
$\big\{ \mu_{y_k,t}\big\}$. 

This completes the proof of the part 1) in Theorem 2.1. 
\vskip.3cm 

\noindent{\bf Step 2:} \, A General Identity. \par 
It remains to analyze the behavior of $\mu$ 
at the end point $y=0$ which shall provide us with the desired boundary 
entropy inequality. We are going to use a general identity
which immediatly follows from the Green formula applied to (2.6).

Let $\theta(t)$ and $\varphi(y)$ be smooth functions not necessarily having 
compact support. 
We multiply the equation (2.6) by $\nabla U(v^\eps) \, \theta \, \varphi$
and integrate over the domain $(y_1, y_2) \times (T_1, T_2)$. 
Integrating by parts and re-ordering the terms, 
we obtain the identity 
$$
E_I^\eps + E_{II}^\eps + E_{III}^\eps = E_{IV}^\eps 
\tag 2.23 
$$
with 
$$
E_I^\eps \equiv - \eps \int_{T_1}^{T_2}\int_{y_1}^{y_2} 
U(v^\eps) \del_t\theta \varphi \, dydt
    + \eps \, \theta(T_2)  \int_{y_1}^{y_2} U(v^\eps(T_2)) \,\varphi \, dy
    - \eps \, \theta(T_1)  \int_{y_1}^{y_2} U(v^\eps(T_1)) \, \varphi \, dy,
\tag 2.24.I
$$
$$
\aligned 
E_{II}^\eps \equiv - \int_{T_1}^{T_2}\int_{y_1}^{y_2} F(v^\eps) 
\theta \del_y\varphi \, dydt
        & + \varphi(y_2) \int_{T_1}^{T_2} \bigl(F(v^\eps(y_2))  
                  - \del_y U(v^\eps)_{|y=y_2} \bigr) \theta  \, dt\\
        & - \varphi(y_1) \int_{T_1}^{T_2} \bigl(F(v^\eps(y_1))
                  - \del_y U(v^\eps)_{|y=y_1} \bigr) \theta  \, dt,
\endaligned 
\tag 2.24.II
$$
$$
\aligned 
E_{III}^\eps \equiv - \int_{T_1}^{T_2}\int_{y_1}^{y_2} U(v^\eps) 
\, \theta \,\del_{yy} \varphi \, dydt
        & + \del_y \varphi(y_2) \int_{T_1}^{T_2} U(v^\eps(y_2)) \,  \theta  \, dt \\ 
        & - \del_y \varphi(y_1) \int_{T_1}^{T_2} U(v^\eps(y_1)) \,  \theta  \, dt,
\endaligned 
\tag 2.24.III
$$ 
and 
$$
E_{IV}^\eps \equiv - \int_{T_1}^{T_2}\int_{y_1}^{y_2} 
       \nabla U(v^\eps)\cdot\bigl(\del_yv^\eps, \del_yv^\eps\bigr) \theta \varphi \, dydt.
\tag 2.24.IV
$$

In case that $\theta \ge 0$ and $\varphi \ge 0$ and 
since $U$ is assumed to convex, one has 
$$
E_{IV}^\eps \le 0, 
\tag 2.25 
$$ 
so we can focus attention on estimating the terms 
$E_I^\eps$, $E_{II}^\eps$ and $E_{III}^\eps$. 
\vskip.3cm 

\noindent{\bf Step 3:} \, Viscous Flux at the Boundary. \par 
We prove here that the viscous flux at the boundary, i.e.~the function 
$\del_y v^\eps(0,t)$, is uniformly bounded in a certain sense
and we determine its weak limit as $\eps \to 0$. 
We use the identity (2.23)-(2.24) with the following choice of parameters:
$$
\supp \theta \subset [T_1,T_2], \quad \supp \varphi \subset [0,1), 
\quad y_1 =0, \quad y_2 =1, \quad (U,F) = (u_j, f_j), \, \, 1 \leq j \leq N. 
$$ 
For $\varphi$ fixed, we obtain
$$
|E_I^\eps | \, \le  \, O(\eps) \, \|\theta\|_{BV}, 
$$
$$
\aligned 
E_{II}^\eps & = - \int_{T_1}^{T_2}\int_0^1 f(v^\eps)\, \theta \, \del_y\varphi \, dydt
        - \varphi(0) \int_{T_1}^{T_2} \bigl(f(u_B^\eps)
                  - \del_y v^\eps(0,.) \bigr) \, \theta  \, dt \\ 
       & =  O(1) \,  \|\theta\|_{L^\infty}
          - \varphi(0) \int_{T_1}^{T_2} \bigl(f(u_B^\eps)
                         - \del_y v^\eps(0,.) \bigr) \, \theta \, dt, 
\endaligned 
$$
and 
$$
\aligned 
E_{III}^\eps & = - \int_{T_1}^{T_2}\int_0^1  v^\eps \theta \del_{yy} \varphi \, dydt
          - \int_{T_1}^{T_2} u_B^\eps \theta \del_y \varphi(0) \, dt \\
        & = O(1) \, \|\theta\|_{L^\infty}. 
\endaligned 
$$ 
Since in this case $E_{IV}^\eps = 0$ and choosing $\varphi$ so that 
$\varphi(0) \ne 0$, it follows 
$$
\bigl|\int_{T_1}^{T_2} \bigl(f(u_B^\eps)
- \del_y v^\eps(0,.) \bigr) \theta  \, dt\bigr|
\le \,  O(1) \, \|\theta\|_{L^\infty} + 
        O(\eps ) \, \|\theta\|_{BV}.
\tag 2.26 
$$

More precisely we can pass to the limit in the identity (2.23) and get
$$
\aligned 
\varphi(0) \, \lim_{\eps \to 0} \int_{T_1}^{T_2} & \bigl(f(u_B)
 - \del_y v^\eps(0,t) \bigr) \theta  \, dt    \\  
= & - \int_{T_1}^{T_2}\int_0^1 <\mu, f> \theta \del_y\varphi \, dydt
    -  \int_{T_1}^{T_2}\int_0^1 <\mu,\id> \theta \del_{yy} \varphi \, dydt 
    -  \del_y \varphi(0) \int_{T_1}^{T_2} u_B \theta  \, dt. 
\endaligned 
$$
On the other hand, it has been observed in Step 1 that (2.20) holds and $<\mu,\id>$
has a trace at $y=0$. Thus one has
$$
\aligned 
  \int_{T_1}^{T_2}\int_0^1 <\mu, f> \theta \del_y\varphi \, dydt
   & + \int_{T_1}^{T_2}\int_0^1 <\mu,\id> \theta \del_{yy} \varphi \, dydt \\ 
&  = 
   \int_{T_1}^{T_2}\int_0^1 C_*(t) \, \theta \del_y\varphi \, dydt
   - \int_{T_1}^{T_2} <\mu_0,\id> \theta \del_y \varphi(0) \, dt \\ 
& =
   - \int_{T_1}^{T_2} C_*(t) \, \theta  \varphi(0) \, dt
   - \int_{T_1}^{T_2} <\mu_0,\id> \theta \del_y \varphi(0) \, dt
\endaligned 
$$ 
and therefore 
$$
\aligned 
\varphi(0) \, \lim_{\eps \to 0} \int_{T_1}^{T_2} & \bigl(f(u_B) 
       - \del_y v^\eps(0,t) \bigr) \theta \, dt\\
        \,  =  \, &  \varphi(0) \int_{T_1}^{T_2} C_*(t) \, \theta  \, dt
       \, + \,\del_y \varphi(0) \int_{T_1}^{T_2} <\mu_0,\id> \theta  \, dt 
      \,- \,\del_y \varphi(0) \int_{T_1}^{T_2} u_B \theta \, dt.
\endaligned 
$$ 
Choosing two test-functions $\varphi$, 
one such that $\varphi(0) =0$ but $\del_y\varphi(0) \ne 0$,
and the other such that
$\varphi(0) \ne 0$ but $\del_y\varphi(0) = 0$, we deduce from the above formula that 
$$
\aligned 
 \lim_{\eps \to 0} \int_{T_1}^{T_2} \bigl(f(u_B) 
- \del_y v^\eps(0,t) \bigr) \theta \, dt
  & = \int_{T_1}^{T_2} C_*(t) \, \theta  \, dt \\ 
\int_{T_1}^{T_2} <\mu_0,\id> \theta \, dt & = \int_{T_1}^{T_2} u_B \theta \, dt.
\endaligned 
\tag 2.27 
$$ 
The first statement in (2.27) is the desired convergence result. 
The second statement is a first step toward proving (2.9). 
\vskip.3cm 

\noindent{\bf Step 4:} \, Boundary Entropy Inequalities (I). \par 
Using (2.27), we are now able to obtain the boundary entropy inequalities. 
We use the identity (2.23)-(2.24) with 
$$
 \theta \ge 0, \quad \supp \theta \subset [T_1,T_2], \quad 
 \varphi \ge 0, \quad \supp \varphi \subset [0,\infty), 
\quad y_1 =0, \quad y_2> 0, 
$$
and $(U,F)$ arbitrary. We obtain
$$
|E_I^\eps | \, \le  \, O(\eps) \, \|\theta\|_{BV}, 
$$
$$
\aligned 
E_{II}^\eps & = - \int_{T_1}^{T_2}\int_0^{y_2} F(v^\eps) 
\,\theta \,\del_y\varphi \, dydt
        - \varphi(0) \int_{T_1}^{T_2} \bigl(F(u_B^\eps)
                  - \del_y U(v^\eps)_{y=0} \bigr) \,\theta  \, dt \\ 
                   & = - \int_{T_1}^{T_2}\int_0^{y_2} F(v^\eps) 
\,\theta \,\del_y\varphi \, dydt
        - \varphi(0) \int_{T_1}^{T_2} \biggl(F(u_B^\eps)
                  - \nabla U(u_B^\eps) \,\del_y v^\eps(0,.) \biggr) \theta  \, dt \\ 
                  & \rightarrow  - \int_{T_1}^{T_2}\int_0^{y_2} <\mu, F> 
\,\theta \,\del_y\varphi \, dydt
        - \varphi(0) \int_{T_1}^{T_2} \biggl(F(u_B)
                  - \nabla U(u_B) \big(f(u_B) - C_*(.)\big) \biggr) \,\theta  \, dt, 
\endaligned 
$$
where we have used (2.27) and the fact that $u_B^\eps \in BV$ converges strongly to $u_B\in BV$,
and 
$$
\aligned 
E_{III}^\eps & = - \int_{T_1}^{T_2}\int_0^{y_2}  U(v^\eps) \,\theta \,\del_{yy} \varphi \, dydt
          - \int_{T_1}^{T_2} U(u_B^\eps) \,\theta \,\del_y \varphi(0) \, dt. 
\endaligned 
$$

Since $E_{IV}^\eps \le  0$  we pass to the limit in (2.23) and get
$$
\aligned 
\varphi(0) \int_{T_1}^{T_2} & \biggl(F(u_B)
                  - \nabla U(u_B) \big(f(u_B) - C_*(t)\big) \biggr) \, \theta  \, dt,  \\ 
\ge 
& - \, \int_{T_1}^{T_2}\int_0^{y_2} 
              \bigl(<\mu_{y,t}, F> \, \del_y\varphi \,  + <\mu_{y,t},U> \, \del_{yy} \varphi\bigr) 
               \,\theta\, dydt \\ 
& + \, \varphi(y_2) \, \int_{T_1}^{T_2} 
              \bigl(<\mu_{y_2,t}, F> - \del_y <\mu_{y,t},U>_{y=y_2}\bigr) \,\theta\, dt \\ 
&  +  \, \del_y \varphi(y_2) \int_{T_1}^{T_2}  <\mu_{y_2,t}, U> \, \theta  \, dt 
  -  \, \del_y \varphi(0) \int_{T_1}^{T_2} U(u_B) \, \theta  \, dt. 
\endaligned 
$$
On one hand, using the test-function $\varphi(y) \equiv 1$, we deduce that
$$
\int_{T_1}^{T_2} \biggl(F(u_B) - - \nabla U(u_B) \big(f(u_B) - C_*(t)\big) \biggr) \theta \, dt
\,  \ge \, \int_{T_1}^{T_2} \big(<\mu, F> + \del_y <\mu, U>_{y=y_2} \big) \, \theta  \, dt 
\tag 2.28
$$
which proves the first inequality in (2.8). 

On the other hand, using the function $\varphi(y) = y$, we obtain 
$$
\aligned 
0 \, \ge & - \, \int_{T_1}^{T_2}\int_0^{y_2} <\mu_{y,t}, F> \,\theta\, dydt 
+ \, y_2 \, \int_{T_1}^{T_2} 
             \bigl(<\mu_{y_2,t}, F> - \del_y <\mu_{y,t},U>_{y=y_2}\bigr) \,\theta\, dt \\ 
 &  +  \, \int_{T_1}^{T_2}  <\mu_{y_2,t}, U> \, \theta  \, dt 
  -  \, \int_{T_1}^{T_2} U(u_B) \, \theta  \, dt, 
\endaligned 
$$
which as $y_2 \to 0$ yields 
$$
\int_{T_1}^{T_2} U(u_B) \theta \, dt \, \ge \, 
\lim_{y \to 0+} \int_{T_1}^{T_2} <\mu_{y,t},U> \theta  \, dt. 
\tag 2.29  
$$ 
In particular, plugging $(U,F) = (u_j,f_j)$, $1 \leq j \leq N$, 
in (2.29), we recover the second statement in (2.27), 
which used together with (2.29) for any fixed, strictly convex entropy $U$ gives:
$$
\aligned 
& \int_{T_1}^{T_2} <\mu_{0,t},U - U(u_B) - \nabla U(u_B) (\id -u_B)> \theta  \, dt \\
& \,= \, \lim_{y \to 0+} \int_{T_1}^{T_2} <\mu_{y,t},U - U(u_B)
 - \nabla U(u_B) (\id -u_B)> \theta  \, dt  \\ 
& \, \le \, 
\int_{T_1}^{T_2} U(u_B) \theta \, dt  \, - 
\int_{T_1}^{T_2} U(u_B) \theta \, dt \,  \\ 
& \, = \, 0. 
\endaligned 
$$ 
But the function $u \to U(u) - U(u_B) - \nabla U(u_B) (u -u_B)$ is positive 
everywhere except at $u_B$ where it achieves its global minimum value. 
It follows that $\mu_{0,t}$ is a Dirac mass concentrated at $u_B$. That proves (2.9). 
\vskip.3cm

\noindent{\bf Step 5:} \, Boundary Entropy Inequalities (II). \par 
We now establish the third inequalities in (2.8). 
We use once more the identity (2.23)-(2.24) with now 
$$
\theta \ge 0, \quad \supp \theta \subset [T_1,T_2], \quad 
 \varphi \ge 0, \quad \supp \varphi \subset [y_1,\infty), 
\quad y_1>0, \quad y_2 =\infty, 
$$
with a function $\varphi$ depending on $\eps$, that is
$$
\varphi^\eps(y,t) \, \equiv \, \tilde \varphi(\eps y, t)
$$
with $\tilde \varphi$ fixed. In that situation one can check that 
$$
\aligned 
E_I^\eps \, & = \, - \int_{T_1}^{T_2}\int_{\eps y_1}^\infty U(u^\eps) 
\,\del_t\theta \,\tilde\varphi \, dxdt\\
  &\,  \rightarrow \, - \int_{T_1}^{T_2}\int_0^\infty 
            <\nu_{x,t}, U> \,\del_t\theta \,\tilde\varphi \, dxdt, 
\endaligned
$$
$$
\aligned 
E_{II}^\eps 
\, & \, = - \int_{T_1}^{T_2}\int_{\eps y_1}^\infty F(u^\eps) \,\theta \,
\del_x\varphi \, dxdt  \\ 
\, & \,  - \tilde\varphi(\eps y_1) \int_{T_1}^{T_2} \bigl(F(v^\eps) 
- \del_yU(v^\eps)_{|y=y_1} \bigr) \theta\, dt \\
\, & \, \rightarrow  - \int_{T_1}^{T_2}\int_0^\infty <\nu_{x,t}, F>
\,\theta \,\del_x\varphi \, dxdt
        - \tilde\varphi(0) \, \int_{T_1}^{T_2} \bigl(<\mu_{y_1,t}, F> 
		  - \del_y <\mu, U>_{|y=y_1}\bigr) \,\theta \, dt, 
\endaligned 
$$
and 
$$
\aligned 
E_{III}^\eps & = - \, \eps \, \int_{T_1}^{T_2}\int_{\eps y_1}^\infty
U(u^\eps) \,\theta \,\del_{xx} \varphi \, dxdt
         \, - \, \del_x\tilde\varphi(\eps y_1) \, \int_{T_1}^{T_2} 
			U(v^\eps)_{|y=y_1} \,\theta  \, dt \\
\, \rightarrow \, 0. 
\endaligned 
$$ 
Since $E_{IV}^\eps \le  0$  and 
$$
- \int_{T_1}^{T_2}\int_0^\infty <\nu_{x,t}, F> \,\theta \,\del_x\tilde\varphi \, dxdt
\, = \, 
\int_{T_1}^{T_2} <\nu_0, F> \,\theta \, dt \, + \, O(1) \, \|\tilde\varphi\|_{L^1}, 
$$
we obtain an inequality of the form 
$$
\tilde\varphi(0) \, \int_{T_1}^{T_2} \bigl(<\mu_{y_1,t}, F> 
- \del_y <\mu, U>_{|y=y_1}\bigr) \,\theta \, dt
\, \ge \, \tilde\varphi(0) \,
\int_{T_1}^{T_2} <\nu_0, F> \,\theta \, dt \, + \, O(1) \, \|\tilde\varphi\|_{L^1}, 
\tag 2.30 
$$
which proves the third inequality in (2.8) by chosing $\tilde\varphi \ge 0$ such that 
$\|\tilde\varphi\|_{L^1} \to 0$ but $\tilde\varphi(0)>0$. 

This complete the proof of Theorem 2.1. 
$\hfill \Box$
\endproclaim

\proclaim{Remark} \rm 
Additional uniform estimates and regularity can be obtained from the identity 
in Step 2 of the proof of Theorem 2.1. 
Let $(U,F)$ be a non-negative entropy pair that is
uniformly convex on $\UU$. 
Use the identity (2.23)-(2.24) with
$$
\theta \equiv 1, \quad T_1 =0, \quad T_2 =T, \quad 
\quad \varphi \equiv 1, \quad y_1 =0, \quad y_2 =\infty.
$$
We assume additonally here that,  
for a fixed state $u_\infty$ and for all $t$, 
$$
u^\eps(x,t) \to u_\infty, \quad u^\eps_x(x,t) \to 0
\qquad \quad \text{ as } x \to \infty. 
$$
The initial data $u_I$ should also decay rapidly at infinity. 
We obtain the following identity
$$
\aligned 
&  \eps \, \int_0^T U(v^\eps(y,T)) \, dy - \eps \, \int_0^\infty U(v^\eps(y,0)) \, dy 
+ \int_0^T F(u_\infty) \, dt  \\ 
- \int_0^T \bigl(F(u^\eps_B) - \nabla U(u^\eps_B) \del_y v^\eps(0,.) \, dt 
&   + \int_0^T\int_0^\infty \nabla^2 U(v^\eps) 
\cdot\big(\del_y v^\eps, \del_y v^\eps\big) \, dy dt 
\, = \, 0. 
\endaligned 
$$
Since the following two terms are uniformly bounded
$$
\bigl| \eps \, \int_0^\infty U(v^\eps(y,0)) \, dy\bigr| \, = \, 
\bigl| \eps \, \int_0^\infty U(u_I^\eps) \, dx \, \bigr| 
\, \le O(1), 
$$
$$
\bigl| \int_0^T \nabla U(u^\eps_B) \del_y v^\eps(0,.) \, dt \big| 
\, \le \, O(1), 
$$
(Cf.(2.26) with $\theta \equiv 1$), we deduce the uniform bounds
$$
\eps \, \int_0^T U(v^\eps(T)) \, dy 
+ \int_0^T\int_0^\infty 
\nabla^2 U(v^\eps) \cdot\big(\del_y v^\eps, \del_y v^\eps\big) \, dy dt 
\, \le \, O(1). 
\tag 2.31 
$$

For every Lipschitz continuous function $g$, 
it follows from (2.31) that the sequence $\del_y g(v^\eps)$ is bounded in $L^2$,
so converges weakly to a limit which is nothing but $\del_y <\mu,g>$:
$$
\del_y g(v^\eps) \, \to \, \del_y <\mu,g>  \qquad \text{ weak 
in } \, L^2(\RR\times\RR_+). 
\tag 2.32 
$$
$\hfill \Box$
\endproclaim

\vskip.3cm

\noindent{\bf 2.2. Finite Difference Schemes.}  \, 
We now extend the above analysis to 
several classes of finite difference schemes 
that are known to be consistent with the entropy inequality
(1.3). Theorem 2.3 below deals with the entropy flux-splittings introduced by 
Chen-LeFloch \cite{\refChenLeFloch}, which also includes as a special case 
the Lax Friedrichs type schemes. We treat the Godunov scheme
in Theorem 2.4.

We are given two mesh parameters $\tau$ and $h$ with 
$\lambda\equiv \tau/h$ kept constant and small enough in order to garantee
the stability of the scheme. 
We define the approximate solutions $u^h(x,t)$ by the scheme
$$
u^h(x,t+\tau) = u^h(x,t)  
-\lambda g\big(u^h(x,t), u^h(x+h,t)\big) + \lambda g\big(u^h(x-h,t), u^h(x,t)\big)
\tag 2.33
$$
and the initial and boundary conditions:
$$
\aligned 
&  u^h(x,t) = u_I(x) \qquad \text { for all } \, t< \tau, \\ 
&  u^h(x,t) = u_B(t) \qquad \text { for all } \, x< h.
\endaligned 
\tag 2.34 
$$
By convention, the functions $u^h$ are right continuous. 
For the Lax-Friedrichs type schemes, the numerical flux $g$ is given by 
$$
g_{\text{Lax}}(v, w) = \frac 1 2 (f(v) + f(w)) - \frac Q \lambda (w-v), 
\tag 2.35 
$$ 
where $Q \in (0,1)$ is called the numerical coefficient of the scheme. 
(Symmetric positive definite matrices $Q$ could also be dealt with.)
For the flux-splitting schemes, $g$ takes the form
$$
g_{\text{split}}(v,w) = f^-(w) + f^+(v),
\tag 2.36 
$$
where $f = f^- + f^+$ is a given entropy flux-splitting for the system (1.1). 
By definition \cite{\refChenLeFloch}, the matrix $\nabla f^\pm$ have real 
eignevalues and a basis of eigenvectors and 
there exists a pair of functions $F_\pm$ such that $(U,F^\pm)$ is an 
entropy pair for the system associated with flux-functions $f^\pm$. 
Observe that (2.35) is a special case of (2.36) as was pointed out
by Chen-LeFloch.

As in the analysis of Section 2.1, we assume a uniform $L^\infty$ bound: 
$$
\|u^h\|_{L^\infty(\RR\times\RR_+)} \le O(1). 
\tag 2.37
$$
We rescale $u^h$ and define the function $v^h:\RR\times\RR_+ \to \UU$ by
$$
v^h(y,t) = u^h(yh,t) \qquad y \ge 0, \, t \ge 0. 
$$
Let $\nu$ and $\mu$ be two Young measures associated 
with $u^h$ and $v^h$, respectively. 

The entropy flux-splitting schemes satisfy  
discrete entropy inequalities of the form
$$
U(u^h(x,t+\tau)) - U(u^h(x,t+\tau))  
+ \lambda \biggl(G(u^h(x,t), u^h(x+h,t)) - G(u^h(x-h,t), u^h(x,t))\biggr)
\le 0, 
\tag 2.38
$$
where $G$ is called the numerical entropy flux. 
With obvious notation, we have 
$$
G_{\text{Lax}}(v, w) = \frac 1 2 (F(v) + F(w)) - \frac Q \lambda (U(w) - U(v))
\tag 2.35bis
$$ 
and 
$$
G_{\text{split}}(v,w) = F^-(w) + F^+(v).
\tag 2.36bis 
$$
Note that (2.38) hold for (2.36)-(2.36bis) provided $u$ takes its value in a 
sufficiently small neighborhood of a given state in $\UU$. 
This is in constrast with the vanishing viscosity method
where no such assumption was necessary.

Theorem 2.1 admits the following extension to the flux-splitting schemes. 
We omit the proof which follows the lines of the one of Theorem 2.1. 

\proclaim{Theorem 2.3} 
Assume that $\UU$ is a small neighborhood of a constant state in $\RR^N$. 
The measure $\mu_{y,t}$ is defined for all $y \ge 0$ and almost every $t$, and is 
constant for $y \in [k, k+1)$ for any integer $k$.  
For all convex entropy pairs $(U,F)$, 
all $y \ge 0$, and in the sense of distributions in $t \in \RR_+$, 
one has 
$$
\aligned
F^+(u_B) + <\mu_{1,t}, F^->  \,
\ge &  \, <\mu_{y,t}, F^+> + <\mu_{y+1,t}, F^->   \\
\ge &  \, <\mu_{y+1,t}, F^+> + <\mu_{y+2,t}, F^-> \\
\ge &  \, <\nu_{0,t}, F>,  
\endaligned 
\tag 2.39 
$$
$$
\mu_{0,t}  = \delta_{u_B(t)} \qquad \text{ for a.e. } \, t>0, 
\tag 2.40
$$
and
$$
\lim_{y \to + \infty} <\mu_{y,t}, F^+> + <\mu_{y+1,t}, F^->  \, \ge \, <\nu_{0,t}, F>.  
\tag 2.41 
$$
$\hfill\Box$
\endproclaim

Consider next the Godunov scheme corresponding to the flux $g$ given by 
$$
g_{\text{Godunov}}(v,w) = f(R(v,w)), 
\tag 2.42 
$$
where we denote by $R(v,w)$ the value at $x/t =0+$ of the solution to the Riemann 
problem with $v$ and $w$ as left and right initial data, respectively. 
The entropy flux is 
$$
G_{\text{Godunov}}(v,w) = F(R(v,w)), 
\tag 2.42bis 
$$
Here it is more convenient to consider the values $R(u^h(x,t), u^h(x+h,t))$ 
and define a function $w^h$
$$
w^h(y,t) = R(u^h(yh,t), u^h(yh+h,t))
\tag 2.43 
$$
for all $y \ge 0$. We denote by $\pi$ a Young measure associated with 
$w^h$ and by $\nu$ a Young measure for $u^h$. It is not difficult to extend Theorem 2.3 as 
follows:

\proclaim{Theorem 2.4} 
The measure $\pi_{y,t}$ is defined for all $y \ge 1/2$ and almost every $t$, 
and is constant in $y$ for $y \in [k-1/2, k+1/2)$ for any integer $k \ge 1$.  
For all convex entropy pairs $(U,F)$, 
all $y \ge 1/2$, and in the sense of distributions in $t \in \RR_+$, 
one has 
$$
\aligned
<\pi_{1/2,t}, F> 
\ge &  \, <\pi_{y,t}, F>  \\
\ge &  \, <\pi_{y+1,t}, F> \\
\ge &  \, <\nu_{0,t}, F>, 
\endaligned 
\tag 2.44
$$
and, at $y = 1/2$ and $y =\infty$, $\pi$ satisfies
$$
<\pi_{1/2,t}, F>  = \lim_{h \to 0}  R(u_B, v^h(1,t)), 
\tag 2.45 
$$
and 
$$
\lim_{y \to \infty} <\pi_{y,t}, F> \, \ge \, <\nu_{0,t}, F>.
\tag 2.46 
$$
$\hfill\Box$
\endproclaim

We conclude this section by giving the main conditions satisfied by the discrete
boundary layer, which will be studied in the rest of this paper. 

Assuming in the results of Theorem 2.3 that $\mu$ is a Dirac mass, 
say $\mu = \delta_v$, 
the {\it discrete boundary layer equation\/} associated with 
the scheme (2.33) takes the form: 
$$
\aligned 
&  g(v(y-1),v(y)) - g(v(y), v(y+1)) = 0 \qquad \text{ for all } \, y \ge 1, \\ 
&   v(y) = u_B, \qquad y \in [0,1), 
\endaligned 
\tag 2.47 
$$
while the {\it discrete boundary entropy inequality\/} is
$$
G(u_B,v_1) \ge F(u_0), 
\tag 2.48 
$$
where $v_1$ plays the role of a parameter. 
Formally, Theorem 2.4 leads to the same equations (2.47)-(2.48) with 
flux and entropy-fluxes given by (2.42).

\subheading {3. Sets of Admissible Boundary Values}

Based on the results in Section 2, 
we introduce in this section several sets which can be used to 
formulate the boundary condition. 
For every method of approximation considered in Section 2, 
we introduce two different sets of admissible boundary values:
\roster 
\item One is based on the entropy inequalities,  
$\EE^{\text{entropy}}(u_B)$ and yields a boundary condition of the form (1.5).
This boundary condition is rigorously satisfied by the limiting function generated
by a sequence of approximate solution. 
as was proven in Section 2.  
For arbitrary systems having few or even just one entropy, 
the set $\EE^{\text{entropy}}(u_B)$ may be too large 
to lead to a well-posed problem; 
\item Another set, $\EE^{\text{layer}}(u_B)$,  
is based on the boundary layer equation, which was obtained formally after 
the analysis in Section 2. 
This set is more adapted to deal with general systems and lead to a 
well-posed problem when the boundary is non characteristic. 
\endroster 

In this section, we study the local structure of those 
sets; under certain assumptions, we can 
prove that the sets $\EE^{\text{layer}}(u_B)$
are manifolds with dimension equal to the number of negative wave speeds
of the system (1.1). This ensures that the initial-boundary 
value problem is well posed if, for instance, the data are constant states 
(boundary Riemann problem) as can be seen by applying the theory in 
\cite{\refLiYu}. 
We recall that (1.1) is assumed to be strictly hyperbolic throughout 
this section and 
we denote by $\lambda_j(u)$ the $N$ real and distinct 
eigenvalues of the matrix $\nabla f(u)$ and by 
$\ell_j(u)$ and $r_j(u)$ corresponding basis of left and right eigenvectors.

\subheading{3.1 Vanishing Viscosity Method}
For the sake of generality, we consider  
$$
\del_t u^\eps + \del_x f(u^\eps) = 
\eps \, \del_x \big(B(u^\eps)\del_x u^\eps\big), 
\qquad \quad x>0, \, t > 0.
\tag 3.1
$$
Theorem 2.1 could be partially extended to this case. 
We assume that the viscosity matrix $B(u)$ depends smoothly upon
its argument $u$ and is positive. We consider entropies $U$ 
that are {\it $B$-convex\/} in the sense that  
$\nabla^2U (u) B(u) > 0$ for all $u$ under consideration. 
The boundary layer equation here takes the form
$$
\del_y f(v) = \del_y\bigl(B(v)\del_y v\bigr) 
\tag 3.2
$$
and the boundary entropy inequalities have the same form (2.16) but now 
$U$ must be $B$-convex.

Following Dubois-LeFloch 
\cite{\refDuboisLeFlochtwo}, we introduce a set based on 
the boundary entropy inequalities.  
{}From now on, the time-dependence may be omitted.

\proclaim{Definition 3.1} Given $u_B \in \UU$, 
the set of admissible boundary values based on the entropy inequalities 
associated with the vanishing viscosity method $(3.1)$ is  
$$
\EE_{\text{viscosity}}^{\text{entropy}}(u_B)
= \big\{ u_0 \in \UU; \, \text{ for all $B$-convex } (U,F), \, 
F(u_B) + \nabla U(u_B) 
\bigl(f(u_0) - f(u_B)\bigr) \ge F(u_0) \big\}. 
\tag 3.3
$$
$\hfill\Box$
\endproclaim

It is obvious that this set may be quite large when the system (1.1) only 
admits few entropies. 
For most systems ($N \ge 3$), this set is
too large to be used to formulate the boundary condition. 
In any case, it is difficult to get information on its local
structure at $u_B$. 
For general systems, the following observation is immediate. 
Fix a state $u_B \in \UU$ and suppose that for some $p$ one has
$$
\lambda_p(u_B) < 0 < \lambda_{p+1}(u_B) 
\tag 3.4
$$
and the basis $r_j(u)$ is a family of eigenvectors for $B(u)$. 
Then the set obtained by formally plugging the expansion  
$$
\aligned 
&  f(u_0) \, \thickapprox \,  f(u_B) + \nabla f(u_B) (u_0 - u_B) 
             + \nabla^2 f(u_B) \cdot\big(u_0 - u_B, u_0 - u_B\big),  \\
&  F(u_0) \, \thickapprox  \,  F(u_B) + \nabla F(u_B) (u_0 - u_B) 
             + \nabla^2 F(u_B) \cdot \big(u_0 - u_B, u_0 - u_B\big) 
\endaligned
\tag 3.5
$$
in the definition of $\EE_{\text{viscosity}}^{\text{entropy}}(u_B)$ 
contains $u_B + $ the span of $r_j(u_B),j=1,...p$ and is contained in a cone 
with vertex $u_B$.
Indeed the inequality under consideration in (3.3) then 
 becomes
$$
\nabla^2 U(u_B) \, \nabla f(u_B) \big(u_0 - u_B, u_0 - u_B\big) \, \le \, 0.
\tag 3.6 
$$
Since $U$ is a convex entropy, the eigenvalues of $\nabla f(u_B)$
satisfy (3.4), and $\nabla^2 U(u_B)$ is a positive definite matrix, 
our claim follows.

\vskip.3cm 

We also consider a second set of admissible boundary values, first
introduced by Gisclon and Serre \cite{\refGisclonSerre}.

\proclaim{Definition 3.2}
Given any $u_B \in \UU$, 
the set of admissible boundary values $\EE_{\text{viscosity}}^{\text{layer}}(u_B)$, 
based on the boundary layer equation 
associated with the vanishing viscosity method 
is 
the set of all $v_\infty \in \UU$ such that the problem 
$$
\aligned 
&  \, B(v)\del_y v = f(v) -f(v_\infty), \\
&  \, v(0) = u_B,\\
&  \lim_{y \to \infty} v(y ) = v_\infty. 
\endaligned
\tag 3.7 
$$
admits a (smooth) solution $v(y) \in \UU$ for $y \ge 0$. $\hfill\Box$
\endproclaim

To study the local structure of  
$\EE_{\text{viscosity}}^{\text{layer}}(u_B)$, we apply the following theorem 
concerning the
existence of invariant manifolds. Cf.~Hartman \cite{\refHartman}
for a proof.

\proclaim{Theorem 3.3}
Consider the differential equation 
$$
\frac{d\xi}{dy}  = E \xi + H(\xi, \xi_0),  \quad \xi(y) \in \RR^N, \, y \in \RR, 
\tag 3.8
$$
where $H: \RR^N\times \RR^N \to \RR^N$ is of class $C^1$ and for each $\xi_0$ 
$$
H(0,\xi_0) = \frac{d H}{d\xi} (0,\xi_0) = 0, 
\tag 3.9
$$
and $E$ is a constant square matrix with $d$ eigenvalues
having negative real part, $e$ eigenvalues having positive real part, and
$N-d-e$ eigenvalues having zero real part.  
For every (small enough) $\xi_0 \in \RR^N$, 
let $\xi_y = \xi (y;\xi_0)$ be the solution of (3.7) with 
the initial condition $\xi (0;\xi_0) = \xi_0$.
Denote by $T_y$ the mapping $\xi_0 \to \xi (y;\xi_0)$.  

There exists a one-to-one mapping of class $C^1$, 
$S: \xi \rightarrow S(\xi) = (w^I, w^{II}, w^{III})$, 
having non-vanishing Jacobian and 
defined on a neighborhood of $\xi= 0 \in \RR^N$ onto a neighborhood of 
$(w^I, w^{II}, w^{III})$ $= (0,0,0)$ $\in \RR^d\times\RR^{N-d-e}\times \RR^e$, 
such that the mapping $S T_y S^{-1}$ takes the simple form 
$$
\aligned 
ST_y S^{-1}: \quad 
              w^I_y     &  = e^{P^I y} w^I_0       + W^I    (y; w^I_0, w^{II}_0, w^{III}_0), \\
               w^{II}_y  &  = e^{P^{II} y} w^{II}_0 + W^{II} (y; w^I_0, w^{II}_0, w^{III}_0), \\
			 	w^{III}_y &  
= e^{P^{III}y} w^{III}_0 + W^{III}(y; w^I_0, w^{II}_0, w^{III}_0), 
\endaligned 
\tag 3.10
$$ 
where $P^I$, $P^{II}$, and $P^{III}$ are constant real-valued matrices 
with all eigenvalues having moduli less than one 
so that the matrix exponentials
$e^{P^I}$, $e^{P^{II}}$, and $e^{P^{III}}$ are well-defined, 
the absolute value of any eigenvalue of $e^{P^I}$ is less than 1, and 
that for $e^{P^{III}}$ is greater than 1, and that for $e^{P^{II}}$ is exactly 1.
Moreover the mapping $W^I$, $W^{II}$, and $W^{III}$ are of class $C^1$ 
and their first order partial derivatives 
with respect to $(w^I_0, w^{II}_0, w^{III}_0)$ 
vanish at $(0,0,0)$.  Moreover one has 
$$
W^I = 0  \quad \text{ and }  \quad W^{II} = 0 \qquad 
\text{ if }   \quad w^I_0 = 0 \quad\text{ and } \quad w^{II}_0 = 0,
\tag 3.11
$$ 
and 
$$
W^{II} = 0  \quad \text{ and }  \quad W^{III} = 0 \qquad 
\text{ if }  \quad w^{II}_0 = 0 \quad\text{ and } \quad w^{III}_0 = 0. 
\tag 3.12
$$ 
$\hfill \Box$
\endproclaim

The condition (3.11) means that the
$e$-dimensional plane $\big\{w^I_0 = 0, w^{II}_0 = 0\big\}$ 
is a locally invariant manifold.
If $S(\xi_0)$ belongs to this plane, then 
$|\xi (y;\xi_0)| \to \infty$ as $y \to \infty$.  
The manifold $\big\{\xi \, / \,  w^I_0 = 0, w^{II}_0 = 0\big\}$ is 
called the {\it unstable manifold\/} of initial data for 
the equation (3.8).

The condition (3.12) means that
the $d$-dimensional plane  $\big\{w^{II}_0 = 0, w^{III}_0 = 0\big\}$ is a locally
invariant manifold.If $S(\xi_0)$ belongs to this plane, then 
$\xi(y;\xi_0)\to 0$ as $y\to\infty$. 
The manifold $\big\{\xi \, / \, w^{II}_0 = 0, w^{III}_0 = 0\big\}$ is called the 
{\it stable manifold.\/}  

Using Theorem 3.3 we prove the following result.

\proclaim{Theorem 3.4} Let $u_B \in \UU$ be given and assume that, for all $u$ 
in a small neighborhood of $u_B$, 
$$
\aligned 
& \text{the basis } r_j(u) \text{ is a family of eigenvectors for } B(u), \\
& \text{the eigenvalues of } B(u), \text{ say }  b_j(u), \text{ are positive,}
\endaligned 
\tag 3.13
$$
and 
$$
\lambda_p(u) < 0 \le \lambda_{p+1}(u)
\tag 3.14
$$
holds for some $p$. 
Then the set $\EE_{\text{viscosity}}^{\text{layer}}(u_B)$ contains the point $u_B$
and, locally nearby $u_B$, contains a manifold with dimension $p$ at least. 
When $0 < \lambda_{p+1}(u_B)$, $\EE_{\text{viscosity}}^{\text{layer}}(u_B)$ 
is a manifold with dimension exactly $p$ and its tangent space at the point $u_B$ 
is spanned by the eigenvectors $r_j(u_B)$, $j=1, 2, \cdots, p$. 
$\hfill \Box$
\endproclaim

Assumption (3.13) holds for the examples considered later in this paper, 
but could easily be relaxed. 
A result similar to our Theorem 3.4 is also proved by Gisclon in \cite{\refGisclontwo},  
by another method. 

\proclaim{Proof of Theorem 3.4} \rm 
The system in (3.6) can be written in the form 
$$
\aligned 
&   \frac{d\tilde v}{dy}
     = B(v_\infty)^{-1} \nabla f(v_\infty)\tilde v + G(\tilde v, v_\infty), \\
&   \tilde v(0) = u_B - v_\infty, \\
&   \tilde v(\infty) = 0,  
\endaligned
\tag 3.15
$$
where 
$\tilde v(y) = v(y) - v_\infty$ 
and the mapping $G(\tilde v, v_\infty)$ satisfies $G(0, v_\infty) = 0, \frac{\del
G}{\del \tilde v} (0, v_\infty) = 0$.   
In view of the assumption (3.13), the two matrices 
$\nabla f(v_\infty)$ and $B(v_\infty)^{-1} \nabla f(v_\infty)$ 
have the same eigenvectors, and so exactly the same
number of positive, zero, and negative eigenvalues.  
Let
$$
\hat\lambda_j(v_\infty) = b_j(v_\infty)^{-1} \, \lambda_j(v_\infty)
$$
be the eigenvalues of $B(v_\infty)^{-1} \nabla f(v_\infty)$.  
Applying Theorem 3.3 with 
$$
\xi(y;\xi_0) = \tilde v(y; u_B - v_\infty), 
$$
we see that there exists a one-to-one $C^1$ mapping $S$, 
defined on a neighborhood of $0\in \RR^N$, 
onto a neighborhood of 
$(w^I, w^{II}, w^{III}) = (0,0,0) \in \RR^p\times \RR^{N-p-1} \times \RR^1$,
such that the manifold 
$$
\EE  \, \equiv \, \big\{\tilde v \, / \,  w^{II}(\tilde v) =0, 
\quad w^{III} (\tilde v) = 0  \big\},
$$ 
which is of dimension $p$, is stable.  For any point $u_B - v_\infty$ taken in
this manifold as an initial data for the differential equation in (3.15), the 
solution $\tilde v(y)$ converges to $0$ as $y \to \infty$,  
which is the third condition required in (3.15).

If $v_\infty$ belongs to this manifold, then (3.15) has a solution
and hence $v_\infty$ solves the boundary layer problem.  
Furthermore the local
structure of the set nearby $u_B$ can be described as follows.

Suppose that $0 < \lambda_{p+1}(u_B)$. The following estimate follows from (3.15): 
$$
\tilde v(y) = \sum_{j=1}^N e^{\hat\lambda_i y} 
            \ell_j(v_\infty) \cdot (u_B - v_\infty) 
             r_j (v_\infty) + 0 (\tilde v(y))^2.
\tag 3.16
$$
For the right handside of (3.16) to go to zero,
we must have
$$
g_j(v_\infty) \equiv \ell_j(v_\infty)\cdot(u_B - v_\infty) = 0,  
\qquad j = p+1, \cdots N. 
\tag 3.17
$$
Keeping $u_B$ fixed, consider the map $g: \UU \to \RR^{N-p}$\
with components $g_j$ given by (3.17). We have 
$$
\frac{dg}{dv_\infty}(u_B) =
- (\ell_{p+1}(u_B), \cdots, \ell_N(u_B)),
\tag 3.18
$$ 
whose rank is $N-p$.  
By the implicit function theorem, (3.17) defines a manifold passing 
through $u_B$ and of dimension $p$. By construction its tangent space at $u_B$ coincides
with the one for the stable manifold $\EE$.  
Therefore, in view of (3.18), the tangent space 
at $u_B$ for $\EE$ is spanned by the $r_j(u_B)$, $j=1, 2, \cdots, p$.
$\hfill \Box$
\endproclaim
\vskip.3cm

A general inclusion can be proven regarding the sets introduced in the 
previous sections. 
It has been first pointed out by Serre \cite{\refSerre}
(cf.~also \cite{\refGisclontwo}) 
that: 

\proclaim{Proposition 3.5} 
The two family of sets introduced in Definitions 3.1 and 3.2 satisfy
the inclusion 
$$
\EE_{\text{viscosity}}^{\text{layer}}(u_B) \subset 
\EE_{\text{viscosity}}^{\text{entropy}}(u_B) 
\tag 3.19  
$$
for all $u_B \in \UU$. $\hfill\Box$
\endproclaim

\proclaim{Proof of Proposition 3.5} \rm 
Let $v_\infty$ be a point in $\EE_{\text{viscosity}}^{\text{layer}}(u_B)$
and denote by $y \to v(y)$ the associated boundary layer function
which satisfies $v(0) = u_B$ and $v(\infty) = v_\infty$. 
Consider the following function of the variable $y>0$:
$$
\Omega(y) \equiv F(v_\infty) - F(v(y)) 
- \nabla U(v(y)) \bigl(f(v_\infty) - f(v(y))\bigr). 
$$
It is easy to see that 
$$
\aligned 
\frac {d\Omega}{dy}(y) 
&  = \nabla^2 U(v(y)) \biggl(f(v_\infty) - f(v(y)), f(v_\infty) - f(v(y))\biggr) \\ 
&  \ge  0
\endaligned 
$$
So the function $\Omega$ is non-decreasing, 
and since $\lim_{y \to \infty} \Omega(y) = 0$, we deduce
that 
$\Omega (y) \le 0$ for all $y$, in particular for $y=0$, that is
$$
F(v_\infty) - F(u_B) - \nabla U(u_B) \bigl(f(v_\infty) - f(u_B)\bigr) \le 0.  
$$
Thus $v_\infty$ belongs to  $\EE_{\text{viscosity}}^{\text{entropy}}(u_B)$. 
$\hfill \Box$
\endproclaim

\noindent{\bf 3.2 Finite Difference Schemes.} 
We now turn to formulations of the boundary condition that are   
based on finite difference approximations. We use the notation in Section 2.2.
We consider a scheme
characterized by its mesh parameters $\tau$ and $h$ with $\lambda = \tau/h$ small enough, 
and by its numerical flux $g(.,.)$ and its family of numerical entropy fluxes $G(.,.)$. 
It is tacitly assumed that the values $u$ remain in a small neighborhood of a given state
and attention is restricted to those entropies $U$ such that the discrete entropy
inequalities (2.38) are satisfied. In fact attention is mostly restricted to 
the Lax-Friedrichs type schemes and the Godunov scheme.

\proclaim{Definition 3.6} Given $u_B \in \UU$, 
the set of admissible boundary values based on the entropy inequalities 
associated with difference scheme is  
$$
\EE_{\text{scheme}}^{\text{entropy}}(u_B)
= \big\{ u_0 \in \UU; \, \text{ there exists } v_1 \text{ s.t. for all convex } (U,F), \, 
G(u_B,v_1) )  \ge F(u_0) \big\}. 
\tag 3.20 
$$ 
$\hfill \Box$
\endproclaim

As for $\EE_{\text{scheme}}^{\text{entropy}}(u_B)$, this set may be too large to 
garantee that the boundary value problem is well posed. 
We also use the obvious notation 
$\EE_{\text{Lax}}^{\text{entropy}}(u_B)$, 
$\EE_{\text{splitting}}^{\text{entropy}}(u_B)$, 
and $\EE_{\text{Godunov}}^{\text{entropy}}(u_B)$. 

For general systems and the diagonalizable splittings, 
i.e.~those such that the vectors $r_j$ form a basis of eigenvectors 
for the matrices $\nabla f^\pm$, we have the following fact.  
Consider a Lax-Friedrichs type scheme or, more generally a
diagonalizable, entropy flux-splitting scheme. Fix a state 
$u_B \in \UU$ and suppose that (3.4) holds for some $p$. 
Then the set obtained by formally linearizing the inequalities 
in the definition of $\EE_{\text{scheme}}^{\text{entropy}}(u_B)$ 
contains $u_B +$the span of $r_j(u_B),j=1,...p$ and is contained in a 
cone with the vertex at $u_B$.To see this we formally plug the
second order expansion 
$$
F^\pm(u_0) \, \thickapprox  \,  F^\pm(u_B) + \nabla F^\pm(u_B) (u_0 - u_B) 
             + \nabla^2 F^\pm(u_B) \big(u_0 - u_B, u_0 - u_B\big) 
\tag 3.21
$$
and obtain the second order version of the inequalities in (3.20): 
$$
\nabla F(u_B) (u_0 - u_B) + \nabla^2 F(u_B) \big(u_0 - u_B, u_0 - u_B\big) 
\, \le \, 
\nabla F^-(u_B) (v_1 - u_B) + \nabla^2 F^-(u_B) \big(v_1 - u_B, v_1 - u_B\big). 
$$
Using the trivial entropies (i.e.~choose for $F$ the components of $f$), 
we get an (second order) expression for $v_1$:
$$
\nabla f^-(u_B) (v_1 - u_B) + \nabla^2 f^-(u_B) \big(v_1 - u_B, v_1 - u_B\big)
\, = \, 
\nabla f(u_B) (u_0 - u_B) + \nabla^2 f(u_B) \big(u_0 - u_B, u_0 - u_B\big), 
$$
which can be used to rewrite the above inequality:
$$
\nabla^2 U(u_B) \nabla f(u_B) \big(u_0 - u_B, u_0 - u_B\big)
\, \le \, 
\nabla^2 U(u_B) \nabla f^-(u_B) \big(v_1 - u_B, v_1 - u_B\big). 
$$
At the first order, $v_1$ is given by 
$$
\nabla f^-(u_B) (v_1 - u_B)
\, = \, 
\nabla f(u_B) (u_0 - u_B) 
$$
so we arrive at the inequality
$$
- \nabla f^+(u_B)^T \nabla f^-(u_B)^{-T} \nabla ^2 U(u_B) 
\nabla f(u_B)\big(u_0 - u_B, u_0 - u_B\big)  \, \le \, 0. 
$$
The desired result follows immediatly since $r_j$ is a basis of eigenvectors 
for the matrices $\nabla f^+$, $\nabla f^-$, and $\nabla f$,
and the function $U$ is convex. 
 
\vskip.3cm 
The second family of sets is now defined.

\proclaim{Definition 3.7}
Given any $u_B \in \UU$, 
the set of admissible boundary values $\EE_{\text{scheme}}^{\text{layer}}(u_B)$, 
based on the boundary layer equation 
associated with the difference scheme  
is 
the set of all $v_\infty \in \UU$ such that the problem 
$$
\aligned 
&  \, g(v(y), v(y+1)) \, = \, f(v_\infty), \\
&  \, v(y) = u_B \quad \text {for } y \in [0,1), \\
&  \lim_{y \to \infty} v(y ) = v_\infty, 
\endaligned
\tag 3.22
$$
admits a (piecewise constant) solution $v(y) \in \UU$ for $y \ge 0$.  
$\hfill \Box$
\endproclaim

To study the local structure of 
$\EE_{\text{scheme}}^{\text{layer}}(u_B)$, we apply the following theorem 
concerning the existence of discrete invariant manifolds. 
(Cf.~Hartman \cite{\refHartman} for a proof.)

\proclaim{Theorem 3.8} 
Let $T : \RR^N \to \RR^N$, $ \xi_0 \to \xi_1$, be a mapping of the form 
$$
\xi_1 = \Gamma \xi_0 +  E(\xi_0), 
\tag 3.23
$$ 
where $E(\xi_0)$ is of class $C^1$ for small $\xi_0$ and
satisfy $E(0)=0$ and ${DE \over D \xi_0} (0) = 0$, and 
the matrix $\Gamma$ is constant, non-singular, and has $d \ge 0$, 
$N-d-e$, $e \ge 0$
eigenvalues of absolute value less than $1$, equal to $1$,
and greater than $1$, respectively.  

There exists a map $S$ of a neighborhood of $\xi_0 = 0$ onto a neighborhood of 
the origin in the space of $(w^I_0,w^{II}_0,w^{III}_0)$
$\in \RR^d\times\RR^{N-d-e}\times \RR^e$ 
such that $S$ is of class $C^1$ with non-vanishing Jacobian
and $STS^{-1}$ takes the simple form 
$$
\aligned 
ST S^{-1}: \quad 
                w^I_1     &= A^I     w^I_0       + W^I (w^I_0, w^{II}_0, w^{III}_0), \\
                w^{II}_1  &= A^{II}  w^{II}_0 + W^{II} (w^I_0, w^{II}_0, w^{III}_0), \\
                w^{III}_1 &= A^{III} w^{III}_0 + W^{III}(w^I_0, w^{II}_0, w^{III}_0), 
\endaligned 
\tag 3.24
$$ 
where $P^I$, $P^{II}$, and $P^{III}$ are $d \times d$, $(N-d-e) \times (N-d-e)$, and
$e \times e$ square matrices with eigenvalues of absolute value less than $1$, equal to
$1$, greater than $1$, respectively,
and the mapping $W^I$, $W^{II}$, and $W^{III}$ are of class $C^1$ 
and their first order partial derivatives with respect to $(w^I_0, w^{II}_0, w^{III}_0)$ 
vanish at $(0,0,0)$. 
Moreover one has 
$$
W^I = 0  \quad \text{ and }  \quad W^{II} = 0 \qquad 
\text{ if }  \quad w^I_0 = 0 \quad\text{ and } \quad w^{II}_0 = 0,
\tag 3.25
$$ 
and 
$$
W^{II} = 0  \quad \text{ and }  \quad W^{III} = 0 \qquad 
\text{ if } \quad w^{II}_0 = 0 \quad\text{ and } \quad w^{III}_0 = 0. 
\tag 3.26
$$  
$\hfill \Box$
\endproclaim

The condition (3.25) means that the plane $v_0 = 0, w_0 = 0$ of
dimension $d$ is locally invariant manifold and if $R(\xi_0)$ belongs to 
this manifold then $T^n\xi_0 \to 0$ as $n \to \infty$.  

The
condition (3.26) means that the plane $u_0 = 0, w_0 = 0$ is a
locally invariant manifold and if $R(\xi_0)$ belongs to this manifold,
$\mid T^n \xi_0\mid \to \infty$ as $n \to
\infty$.  

Using this theorem we shall prove: 

\proclaim{Theorem 3.9} Consider a Lax-Friedrichs type scheme. 
Let $u_B \in \UU$ be given and assume that $(3.14)$ holds for some $p$. 
Then 
the set $\EE_{\text{Lax}}^{\text{layer}}(u_B)$ contains the point $u_B$
and, locally nearby $u_B$, contains a manifold with dimension $p$. 
When $0 < \lambda_{p+1}(u_B)$, $\EE_{\text{Lax}}^{\text{layer}}(u_B)$ 
is a manifold with dimension exactly $p$ and its tangent space at the point $u_B$ 
is spanned by the eigenvectors $r_j(u_B)$, $j=1, 2, \cdots, p$. 
$\hfill \Box$
\endproclaim

\proclaim{Proof of Theorem 3.9} \rm We search for all $v_\infty$ that solve the 
problem: 
$$
\aligned 
& H(v(y), v(y+1), v_\infty) \, = \, 0  \\
& v(0) = 0, \\
& v(\infty) = v_\infty
\endaligned 
\tag 3.27 
$$
with 
$$
H(v(y), v(y+1), v_\infty) \equiv v(y+1) - v(y) 
- \frac {\lambda}{2 Q} \bigl( f(v(y)) + f(v(y+1)) -  2 f(v_\infty)\bigr). 
\tag 3.28 
$$ 
Using the notation $H=H(v,w,v_\infty)$, we compute
$$
\aligned 
& \frac {\del H}{\del v} (v,w,v_\infty)  = Id + \frac {\lambda}{2 Q} \nabla f(v), \\
& \frac {\del H}{\del w} (v,w,v_\infty)  = Id - \frac {\lambda}{2 Q} \nabla f(w).  \\
\endaligned
\tag 3.29 
$$ 
For $\lambda/(2Q)$ small enough, the matrix $\del H/\del w$ is invertible
and its inverse is uniformly bounded w.r.t the variables $v$, $w$, and 
$v_\infty$. 
By the global implicit function theorem (see J.T. Schwartz \cite{\refSchwartz}) 
the system (3.27) can be solved for $v(y+1)$. 
So there exists a smooth mapping $K(v(y), v_\infty)$ such that
$$
v(y+1) \, = \, K(v(y), v_\infty) 
\tag 3.30 
$$ 
and  $K(v_\infty, v_\infty) =0$. Moreover one has 
$$
{\del K\over \del v}(v(y), v_\infty)
 = \bigl(Id - \frac {\lambda}{2 Q} \nabla f(v(y+1))\bigr)^{-1} 
                       \bigl(Id + \frac {\lambda}{2 Q} \nabla f(v(y))\bigr). 
\tag 3.31 
$$

The system (3.30) can be linearized around $v_\infty$: 
$$
\aligned
v(y+1) = & -\big(Id - \frac {\lambda}{2 Q} {\del f\over \del u} (v_\infty)\big)^{-1} 
\big(Id + \frac {\lambda}{2 Q} {\del f\over \del u} (v_\infty)\big)v(y) \\ 
& + K(v(y), v_\infty) 
+ \big(Id - \frac {\lambda}{2 Q} {\del f\over \del u} (v_\infty)\big)^{-1} 
\big(Id + \frac {\lambda}{2 Q} {\del f\over \del u} (v_\infty)\big)v(y).
\endaligned 
$$ 
Set $v^*(y+1) = v(y+1) -v_\infty$. The system can be written as 
$$
\aligned
v^*(y+1)= & -\big(Id - \frac {\lambda}{2 Q} {\del f\over \del u} (v_\infty)\big)^{-1} 
\big(Id + \frac {\lambda}{2 Q} {\del f\over \del u} (v_\infty)\big)v^*(y) \\ 
& + G(v^*(y) + v_\infty, v_\infty) 
+ \big(Id - \frac {\lambda}{2 Q} {\del f\over \del u} (v_\infty)\big)^{-1} 
\big(Id + \frac {\lambda}{2 Q} {\del f\over \del u} (v_\infty)\big)v^*(y).
\endaligned
$$
In other words 
$$
v^*(y+1)= A(v_\infty) v^*(y) + K^*(v^*(y), v_\infty),
\tag 3.32
$$ 
where
$$
A(v_\infty) \equiv  \big(Id - \frac {\lambda}{2 Q} {\nabla f} (v_\infty) \big)^{-1}
\big(Id + \frac {\lambda}{2 Q} {\del f\over \del u}(v_\infty) \big)
\tag 3.33a 
$$ 
and 
$$
K^* \ \text{ and }  {\del K^*\over \del v^*(y)}  \text{ vanish } 
\text{ at } v^*(y) = 0.
\tag 3.33b 
$$ 
\vskip.3cm

We observe that 
$$
\text{ the eigenvalues of the matrix $A(v_\infty)$ are 
${1+\frac{\lambda} {2 Q} \lambda_i(v_\infty)\over1-\frac{\lambda} {2 Q}\lambda_i(v_\infty)}$, } 
\tag 3.34 
$$ 
where 
$\lambda_i(v_\infty)$ are the eigenvalues of $\nabla f(v_\infty)$.  

Namely (3.34) follows from the fact that the following two statements are equivalent~: 
\roster 
\item  $a$ is an eigenvalue of $A(v_\infty)$; 
\item there exists $r \neq 0$  such that $A(v_\infty) r = ar$. 
\endroster

Using the expression (3.33a) of $A(v_\infty)$ and simplifying
the resulting equation, we get
$$
{\nabla f} (v_\infty) r = {2Q(a-1)\over \lambda(1+a)} r.
$$ 
So $a$ is an eigenvalue of $A(v_\infty)$ if and only if 
${2Q(a-1)\over \lambda (1+a)}$ is an eigenvalue of ${\del f}{\del u}
(v_\infty)$ with right eigenvector $r$; so 
$$
{2Q(a-1)\over\lambda(a+1)} = \lambda_i(v_\infty)
\tag 3.35
$$ 
for some $i$ with left eigenvector $\ell_i(v_\infty)$ and right
eigenvector $r_i(v_\infty)$.  
Solving (3.35) for $a$ we get $i$th eigenvalue of $A(v_\infty)$
$$
a_i = {1+\frac{\lambda} {2 Q} \lambda_i(v_\infty)\over
1-\frac{\lambda} {2 Q}\lambda \lambda_i(v_\infty)}. 
\tag 3.36
$$ 
Let $T$ be a matrix which diagonalize ${\nabla f}(v_\infty)$.  
Then the same matrix diagonalize $A(v_\infty)$: 
$$
T A T^{-1} = \diag (a_1, a_2, \cdots a_n).
$$ 
Set $w^*(y+1) = Tv^*(y+1)$, we get
$$
w^*(y+1)= \left(
\aligned
&a_1 \; \\
& \;\;\;\;\; a_2 \;\;\;\; 0 \\
&\;\;\;\;\; \ddots \\ 
&0 \;\;\;\;\;\;\;\;\; a_n 
\endaligned
\right) w^*(y) + G^*(T^{-1} w^*(y), v_\infty)
$$
where $G^*$ and ${\del G^*\over \del w^*(y)}$  are zero at
$w^*(y) = 0$. 

Note that 
$$
a_1 < a_2 < \cdots a_p < 1 \le a_{p+1}  < \cdots < a_n. 
\tag 3.37
$$ 
and 
$$
a_{p+1} = 1 \Leftrightarrow \lambda_{p+1} (v_\infty) = 0.
$$ 
Since all the hypothesis of Theorem 3.8 are satisfied, 
there exists a $p$-dimensional invariant manifold defined near $0$ such that, 
if the data $v^*_0$ belongs to this manifold, then $w^*(y+1) \to 0$ as $y
\to \infty$. 
In fact in terms of the original variable $v(y+1)$,
we have the expansion 
$$
v(y+1) - v_\infty = {\displaystyle{\sum_{j=1}^N}} a_j^y <\ell_j(u), v_b -
v_\infty> r_j(v_\infty) + 0 (\mid v(y+1) - v_\infty\mid)^2. 
\tag 3.38
$$ 

In order for this to go to zero, as $y\to 0$  we must have
$$
< \ell_j(v_\infty), u_B - v_\infty>  = 0, \qquad  j = p+1, \cdots N.
\tag 3.39
$$
This for fixed $u_B$ defines a map from $R^N \to R^{N-p}$ and
whose Jacobian at $u_B = v_\infty$ is the matrix whose $N-p$ rows are
$\ell_j(v_\infty)$.  
Since $\ell_j(v_\infty)$ are linearly independent by
implicit function theorem we deduce that (3.39) defines a $p$
dimensional manifold passing through $u_B$ and if $v_\infty$ is in this
manifold then there exist a solution to (3.29) whose local
structure is given by (3.39). 
$\hfill \Box$
\endproclaim

The following general inclusion can be proven:

\proclaim{Proposition 3.10} 
The two family of sets introduced in Definitions 3.6 and 3.7 satisfy, for
all $u_B \in \UU$, 
$$
\EE_{\text{scheme}}^{\text{layer}}(u_B) \subset \EE_{\text{scheme}}^{\text{entropy}}(u_B).
\tag 3.40 
$$
$\hfill \Box$
\endproclaim

\proclaim{Proof of Proposition 3.10} \rm 
We consider as before a difference scheme that satisfies discrete entropy inequalities. 
For every $v_\infty$ in the set $\EE_{\text{scheme}}^{\text{layer}}(u_B$, there exists
a corresponding boundary layer profile $v(y)$, solution of 
$$
g(v(y), v(y+1)) \, = \, f(v_\infty).
$$
The function $v(y)$ is actually a stationnary solution to the scheme since
$$
v(y) - v(y)  + \lambda \big(g(v(y), v(y+1)) - g(v(y-1), v(y))\big) =0. 
$$ 
Therefore for every convex entropy pair $(U,F)$, 
it satifies the entropy inequality
$$
U(v(y)) - U(v(y))  + \lambda \big(G(v(y), v(y+1)) - G(v(y-1), v(y))\big) \, \le \, 0, 
$$
which is nothing but 
$$
G(v(y), v(y+1)) - G(v(y-1), v(y) \, \le \, 0 
$$
Since $\lim_{y \to \infty} v(y ) = v_\infty$, we get
$$
G(v(y), v(y+1)) \, \ge \, F(v_\infty)
$$
and so with $y=0$, since $v(y) = u_B \quad \text {for } y \in [0,1)$, 
$$
G(u_B,v_1) )  \ge F(u_0) 
$$
with $v_1 = v(1)$. That establishes that $v_\infty$ belongs to the set 
$\EE_{\text{scheme}}^{\text{entropy}}(u_B)$. 
$\hfill \Box$
\endproclaim 

\vskip.3cm

Finally we treat the Godunov scheme. The sets 
$\EE_{\text{Godunov}}^{\text{layer}}(u_B)$ 
and $\EE_{\text{Godunov}}^{\text{entropy}}(u_B)$
are defined by Definitions 3.6 and 3.7. We now prove:

\proclaim{Theorem 3.11} Consider the Godunov scheme and let $u_B \in \UU$ be given.
We have
$$
\EE_{\text{Godunov}}^{\text{layer}}(u_B) \, 
= \, \EE_{\text{Godunov}}^{\text{entropy}}(u_B). 
\tag 3.41 
$$
This set can also be described as the set 
$$
\EE^{\text{Riemann}} (u_B) \, = \, \big\{R(u_B, w) \, / \, w \in \UU\big\},
$$
where $R(u_B, w)$ denotes the value at $x/t =0+$ of the solution of the Riemann
problem with data $u_B$ and $w$ on the left and right, respectively. 
Moreover when $(3.4)$ holds for some $p$, 
the set above contains the point $u_B$ and, locally nearby $u_B$, is
a manifold with dimension $p$ and with tangent space at the point $u_B$ 
spanned by the eigenvectors $r_j(u_B)$, $j=1, 2, \cdots, p$. 
$\hfill \Box$
\endproclaim

Observe that the Godunov scheme does not produce any boundary layer, in the 
sense that 
the layer contains no interior point.

\proclaim{Proof of Theorem 3.11} \rm We recall that the set 
$\EE_{\text{Godunov}}^{\text{layer}}(u_B) $
is defined by the equation
$$
\aligned 
& f(u_B) \, = \, f(R(v(y), v(y+1)),  \\
& v(y) = u_B \quad \text{ for all } \, y \in [0,1), \\
& \lim_{y \to \infty} v(y) = v_\infty, 
\endaligned 
\tag 3.42 
$$
while the set  $\EE_{\text{Godunov}}^{\text{entropy}}(u_B)$
is defined by the inequalities
$$
F(R(u_B, v_1)) \ge F(u_0)  \qquad \text { for all convex pair } \, (U,F) 
\tag 3.43
$$
and for some $v_1 \in \UU$. So it is not hard to see from the definition that
$$
\EE^{\text{Riemann}}(u_B) \subset \EE_{\text{Godunov}}^{\text{layer}}(u_B).
$$

On the other hand the inclusion
$$
\EE_{\text{Godunov}}^{\text{layer}}(u_B) \subset 
\EE_{\text{Godunov}}^{\text{entropy}}(u_B)
$$
also holds in view of Proposition 3.10. 

It remains to show that 
$$
\EE_{\text{Godunov}}^{\text{entropy}}(u_B)
\subset 
\EE^{\text{Riemann}}(u_B). 
$$ 
Consider a pair $(u_0, v_1)$ that solves (3.43). Then we need show that 
there exists $w$ such that
$$
R(u_B,w) = u_0. 
\tag 3.44 
$$
Using the trivial entropies, we get 
$$
f(R(u_B,v_1)) = f(u_0)
$$
which, combined with the inequality (3.43), shows that 
the pair of states $(R(u_B,v_1), u_0)$ is an entropy satisfying, stationary
shock wave. On the other hand the Riemann problem with left and right initial
data $u_B$ and $R(u_B,v_1)$, respectively, contains only waves with non-positive
speeds. Therefore the Riemann solution, with $u_B$ as a left state
and $u_0$ as a right state, only contains waves with non-positive speeds. 
This function takes the value $u_0$ in the whole half-interval $x/t >0$ and thus
$R(u_B, u_0) = u_0$, which proves (3.44) with $w=u_0$. 
$\hfill \Box$
\endproclaim

\subheading{4. Selected Examples }

In this section, we consider first the convex scalar conservation laws
and establish that all the sets introduced in Section 2 are essentially
the same. 
Some remarks are then given for the linear hyperbolic systems. 
Next we return to the scalar equation and treat a non-convex flux
function, showing again that the sets are the same with the exception
of the set based on the boundary layer equations. Finally we 
treat the elastodynamics system and isentropic Euler system.  
\vskip.3cm 

\noindent{\bf 4.1. Scalar Conservation laws: Convex Fluxes.} 
We consider a scalar conservation law with strictly convex flux,
i.e. $f^{''}(u)>0$
and analyze the boundary layer equation. Let $u_*$ be the unique point
such that $f'(u_*) = 0$.  To the state $u_B$, when $u_B \ne u_*$, we associate
the solution $u_B^* \ne u_B$ of the equation $f(u_B^*) = f(u_B)$.  
 
We state here a theorem which says that some of the  sets introduced in
Section 3 coincide in this case. We also recover the formulation of the 
boundary condition discovered by Bardos-Leroux-Nedelec 
\cite{\refBardosLerouxNedelec} and Leroux \cite{\refLeroux}.

\proclaim{Theorem 4.1} Consider a scalar conservation laws with convex flux. \par 
1) \, For any $u_B \in \UU\equiv \RR$, the sets of admissible boundary values 
$\EE_{\text{viscosity}}^{\text{entropy}}(u_B)$, 
$\EE_{\text{Godunov}}^{\text{layer}}(u_B)$, and 
$\EE_{\text{Godunov}}^{\text{entropy}}(u_B)$,
coincide with 
$$
\EE^{\text{Riemann}}(u_B) \,  = \, \cases
\alignedat2 
& \left(-\infty, u_B^*\right]  \cup \big\{u_B\big\} \quad
 && \text{ if }  \ u_B > u_*,   \\ 
& \left(-\infty, u_* \right] \quad
&&  \text{ if } \ u_B \leq u_*.
\endalignedat
\endcases
$$
and 
$$\EE^{\text{layer}}_{\text{viscosity}} \ (u_B) =
\EE_{\text{Riemann}}(u_B) \setminus \big\{u_B^*\big\}$$ 
2) \, Given $u_B \in \UU \equiv [-M, M]$ for a fixed value of $M>0$, we set 
$\|f'\|_\infty = \sup_{ w \in [-8M, 8M]} |f'(w)|$ 
and consider a Lax-Friedrichs type scheme with coefficient $\lambda$ and $Q$ satisfying 
$\|f'\|_\infty \lambda /Q \le 1$, then 
$$
\aligned 
& \EE_{\text{Lax}}^{\text{layer}}(u_B) \cap [-M, M] \, 
=  \EE^{\text{Riemann}}(u_B) \cap [-M, M] \setminus \big\{u_B^*\big\} \\
& \EE_{\text{Lax}}^{\text{entropy}}(u_B) \cap [-M, M] 
= \, \EE^{\text{Riemann}}(u_B) \cap [-M, M] 
\endaligned
$$
$\hfill \Box$
\endproclaim

\noindent{\bf 4.2 Linear Hyperbolic Systems.} \, 

It is not hard to prove that for a linear and strictly hyperbolic system, 
the sets defined in Section 3 are all equivalent when boundary is not 
characterestic. We only consider here the case of the discrete boundary 
layer based on the Lax-Friedrichs scheme.

We also focus attention in this section to establish that 
the restriction (3.12) on the viscosity matrix is essential to our purpose
here, 
as was observed in another context by Majda-Pego \cite{\refMajdaPego}
in their study of traveling wave solutions to (2.1). 
The following example shows a situation where the viscosity matrix
is a positive diagonal matrix, and does not satisfy (3.12), while 
the formulation may lead to a ``wrong'' boundary condition.  

We consider the linear system 
$$
\del_t u + 
\left
(\aligned
&-5\;\;\;\;\;\;5 \\
&-3\;\;\;\;\;\;3 
\endaligned
\right)  \del_x u = \eps \, 
\left(
\aligned
& 5 \;\;\;\;\;\; 0 \\
& 0 \;\;\;\;\;\; 1
\endaligned 
\right) \del_{xx} u.
\tag 4.1
$$
According to our earlier analysis, the boundary layer equation is
$$
\del_{yy} v(y) =
\left
(\aligned
&1/5\;\;\;\;\;\;0 \\
&0\;\;\;\;\;\;\;\;1 
\endaligned
\right)
\left(
\aligned
& -5 \;\;\;\;\;\; 5 \\
& -3 \;\;\;\;\;\; 3
\endaligned 
\right) \del_y v(y), 
$$
i.e. 
$$
\del_{yy} v(y)=  
\left
(\aligned
&-1\;\;\;\;\;\;1 \\
&-3\;\;\;\;\;\;3 
\endaligned
\right) \del_y v(y). 
$$
Integrating this equation once and using $v(+\infty)= v_\infty$, we get
$$
\del_y v(y) = 
\left
(\aligned
&-1\;\;\;\;\;\;1 \\
&-3\;\;\;\;\;\;3 
\endaligned
\right)(v-v_\infty).
\tag 4.2
$$ 
Now the eigenvalues of 
$\left(\aligned 
& -5\;\;\;\;\;5 \\
& -3\;\;\;\;\;3 
\endaligned
\right)$ are $\lambda_1 = -2$ and $\lambda_2 = 0$. 
On the other hand, the eigenvalues of 
$\left(\aligned 
& -1\;\;\;\;\;1 \\
& -3\;\;\;\;\;3 
\endaligned
\right)$ are $\mu_1 = 0$ and $\mu_2 = 2$.  
The solution of (4.2) with the
initial condition $v(0) = v_B-v_\infty$ is 
$$
v(y)-v_\infty = <\bar\ell_1, v_B - v_\infty> \bar r_1
   + <\bar\ell_2, v_B - v_\infty> \bar r_2{e}^{2y},
$$ 
where
$$
\bar\ell_1 = \left(\frac{-3}{2},\frac{1}{2}\right),
\quad \bar\ell_2 = \left(\frac{1}{\sqrt 2} \;\frac{-1}{\sqrt 2}\right), 
\quad \bar r_1 =
\left(\aligned
&1/\sqrt 2 \\
&1/\sqrt 2\endaligned\right), \quad
\bar r_2 = \left(\aligned
& 1/2 \\ 
& 3/2\endaligned\right).
$$ 
In order for $v(y) \to v_\infty$ as $y \to \infty$, we must have
$<\bar\ell_1, v_B -v_\infty> = 0$ and $<\bar\ell_2, v_B -v_\infty>=0$
which means that $v_\infty = v_B$.
This requires that we prescribe $u$ at the boundary.This is
wrong boundary condition for the hyperbolic system 
$$
\del_t u +
\left(
\aligned 
& -5 \;\;\;\;5 \\
& -3 \;\;\;\;3\endaligned\right) \del_x u = 0
$$
because none of the characterestics are entering.

Let us now consider the numerical boundary layer 
for a general linear and strictly hyperbolic
system. 
Set $f(u)  = Au$, where $A$ is a constant matrix.
The boundary layer equation becomes 
$$
\frac{\lambda}{2} A \, v(y+1) + \frac{\lambda}{2} A \, v(y) 
- \frac 1 2 \big(v(y+1) - v(y)\big) = \lambda \, A \, v_\infty, 
\tag 4.3
$$ 
$$
v(0) = v_B, \qquad v(\infty) = v_\infty.
$$ 
For a given $u_B$, we search for the set of states $v_\infty$ 
for which this problem has a solution. 
Set $v^\infty(y) = v(y+1) - v_\infty$.  The
first equation in (4.3) becomes 
$$
(\lambda A - I) v^\infty(y) = - (\lambda A +I) v^\infty(y-1). 
\tag 4.4
$$ 
Let $\ell_j$ and $r_j$ be the left- and right- eigenvectors for
$A$ associated 
with the eigenvalues $\lambda_j$. Set $C^j(y) = <\ell_j, v(y+1)>$.
{}From (4.4) we get
$$
(1-\lambda \, \lambda_j) C(y)^j = (1 + \lambda \, \lambda_j)
C^j(y-1)
$$ 
or
$$
C^j(y)   = \left(\frac{1+\lambda \, \lambda_j}{1 - \lambda
\, \lambda_j}\right) C^j(y-1) 
$$
with
$$
C_0^j = <\ell_j, v_B - v_\infty>.
$$
Integrating this, we get 
$$
C^j(y) = <\ell_j, v_B- v_\infty > \left(\frac{1 + \lambda \lambda_j}{1-\lambda
\lambda_j}\right)^y
$$ 
or
$$
v(y+1) -  v_\infty = v^\infty(y) = {\displaystyle{\sum_{j=1}^n}}
\left(\frac{1+\lambda \lambda_j}{1-\lambda \lambda_j}\right)^y<\ell_j, u_B - 
v_\infty)
r_j.
$$ 

For $v(y+1) \to v_\infty$, we need $<\ell_j, v_B - v_\infty> = 0, j = p+1,
\cdots n$ because $\lambda_1 < \lambda_2 < \cdots \lambda_p < 0 \le 
\lambda_{p+1} <
\cdots \lambda_n$.  This gives correct boundary condition when the 
eigenvalues are not zero; i.e. to prescribe 
$$
<\ell_j, u> \qquad \text{ for } \,  j = p+1, \cdots, N.
$$

\noindent{\bf 4.3 Scalar Conservation Laws: Non-Convex Fluxes.} 

We return to scalar conservation laws but now with non-convex fluxes. 
For definiteness we treat the case of the cubic flux given by 
$$
f(u) = \frac{1}{2} (u^3 - 3u), 
\tag 4.5
$$
which has one minima and one maxima; indeed
$$
f(1) = -1, \,\, f'(1) = 0, \,\, f^{''} (1) = 3, \,\, 
f(-1) = 1, \,\, f'(-1) = 0, \,\, f^{''}(-1) = -3.
$$
For a given $u_B \in \RR$ and the function $f$ given by (4.5), 
we shall need the solution of the equation 
$$
f(u) = f(u_B), u \neq u_B. 
\tag 4.6
$$
If
$u_B < -2$ or $u_B >2$, there is no solution for (4.6).  
If  $u_B \in (-2, -1) \cup (1,2)$, then (4.6) has exactly two solutions.
In this case we denote by $u_B^{\ell}$ and $u_B^s$ the largest
and smallest solutions of (4.6), respectively.  If  $u_B = -2, -1, 1$, or $2$, 
then (4.6) has exactly one solution; namely $1,2,-2$, and $-1$, 
respectively.

For the formulation of the results in this subsection, it will be convenient 
to introduce the following set, which is either the empty set or contains a 
single element: 
$$
E(u_B) \, = \, \cases
\alignedat2 
&\emptyset, \quad
&&\text{if} \ u_B \in (-\infty, -2) \cup [-1, 1] \cup (2, \infty) \\
&\big\{1\big\} \quad 
&&\text{if} \ u_B = -2 \\
& \big\{u_B^s\big\}, \quad
&& \text{if} \ -2 < u_B < -1 \\
& \big\{u_B^{\ell}\big\},  \quad 
&& \text{if} \ 1 < u_B < 2. \\
&\big\{-1\big\} \quad 
&&\text{if} \ u_B = 2. \\
\endalignedat 
\endcases
\tag 4.7 
$$

\proclaim{Theorem 4.2}  Consider the scalar conservation law with the
non-convex flux $(4.5)$. 
\vskip.3cm 
\noindent{1)}  For any $u_B \in \UU = R$, the set of admissible
boundary values $\EE^{\text{entropy}}_{\text{viscosity}}(u_B),
\EE^{\text{layer}}_{\text{Godunov}} (u_B)$, and
$\EE_{\text{Godunov}}^{\text{entropy}} (u_B)$  coincide with 
$$
\EE^{\text{Riemann}} (u_B) = \, \cases 
\alignedat2 
& \big\{u_B\big\}, \quad 
&&\text{if} \ u_B < -2 \\
&\big\{-2, 1\big\}, \quad 
&&\text{if} \ u_B = -2 \\
& [u_B^s, 1] \cup \big\{u_B\big\}, \quad
&& \text {if}  -2 < u_B < -1 \\
& [-1, 1] \, \quad 
&& \text{if} \ -1 \le u_B \le 1 \\
& [-1, u_B^{\ell}] \cup \big\{u_B\big\}, \quad
&& \text{if} \ 1 < u_B < 2 \\
& \big\{u_B\big\}, \quad
&& \text{if} \ u_B >  2 \\
& \big\{2, -1\big\}, \quad 
&&\text{if} \ u_B = 2 
\endalignedat
\endcases
$$ 
and 
$$
\EE^{\text{layer}}_{\text{viscosity}} (u_B) =
\EE^{\text{Riemann}} (u_B) - E(u_B).
$$
2) Given any state $u_B \in \UU = [-M, M]$ for a fixed value $M>2$, we set
$\parallel f' \parallel_\infty = {\sup_{w \in [-8M, 8M]}} \mid
f'(w)\mid$ and consider a Lax-Friedrichs type scheme with
coefficient $\lambda$ and $Q$ satisfying $\parallel  f'
\parallel_\infty \frac{\lambda}{Q} \le 1$. Then 
$$
\aligned
& \EE_{\text{Lax}}^{\text{layer}} (u_B) \cap [-M, M] =
\EE^{\text{Riemann}}(u_B) \cap [-M, M] \setminus E(u_B), \\ 
& \EE_{\text{Lax}}^{\text{entropy}}(u_B) \cap [-M, M] =
 \EE^{\text{Riemann}} (u_B) \cap [-M, M].
\endaligned
$$ 
$\hfill\Box$
\endproclaim  

The proof of this is straightforward and is omitted.

\noindent{\bf 4.4 Nonlinear Elastodynamics.} \, 
The system considered now arises in the modeling of elastic materials 
\cite{\refCourantFriedrichs}: 
$$
\aligned
\del_t v -\del_x u =0, \\
\del_t u - \del_x \sigma(v) =0. 
\endaligned
\tag 4.8 
$$
It describes the evolution of a nonlinear material 
with deformation gradient $v$
and velocity $u$. The stress function $\sigma$ is assumed to be smooth enough 
and satisfy the following conditions:
$$
\sigma'(v) > 0, \qquad v \, \sigma''(v) > 0.
\tag 4.9
$$

Let us discuss the vanishing viscosity approximation
for the viscosity matrix  $B(u) = I$. 
The boundary layer problem to be studied here is 
$$
\aligned
- \, \del_y u \, = \, \del^2_y v, \\
- \, \del_y \sigma(v) \, = \, \del^2_y u, \\
v(0) = v_B, \qquad v(\infty) = v_\infty, \\
u (0) = u_B, \qquad u(\infty) = u_\infty 
\endaligned
\tag 4.10 
$$
We need determine the set of $(v_\infty, u_\infty)$ for which (4.10)
has a solution.  Integrating once the ODE'S and using the boundary
condition at infinity, we get
$$
\del_y v =  u_\infty - u, \qquad  u_y = \sigma(v_\infty) -\sigma(v). 
\tag 4.11
$$
Cross multiplying the equations and integrating, we get
$$
\frac{(u-u_\infty)^2}{2} = \int_{v_\infty}^v (\sigma(s) - \sigma(v_\infty)) \, ds,
$$ 
so 
$$
(u - u_\infty) = \pm \left(\int_{v_\infty}^v
2 \, (\sigma(s)-\sigma (v_\infty)) \, ds\right)^{1/2}. 
\tag 4.12 
$$ 
Note that $\int_{v_\infty}^v (\sigma (s) - \sigma (v_\infty)) ds \ge
0$ because of the condition $\sigma' (v) > 0$.  From (4.12) it
follows that 
$$
v(y) = v_\infty \quad \Leftrightarrow \quad u(y) = u_\infty.
$$
Since we are interested in a solution
connecting  $(v_B, u_B)$ at $y = 0$ to $(v_\infty, u_\infty)$ at $y =
\infty$, we get from (4.11) that either 
$$
\aligned 
&  v_B < v (y) < v_\infty \quad \text{and}  \quad u_B < u(y) < u_\infty \\ 
& \text{or} \\ 
&  v_B > v (y) > v_\infty \quad \text{and} \quad u_B > u(y) > u_\infty.
\endaligned
\tag 4.13 
$$ 
This determines the sign in (4.12):
$$
u = \cases 
\alignedat2 
& u_\infty + \left(\int^v_{v_\infty} 2 \, (\sigma(s)-\sigma(\bar
v)) \, ds\right)^{1/2} \quad && \text{if} \, v > v_\infty \\
& u_\infty - \left(\int_{v_\infty}^v 2 \, (\sigma(s) -\sigma(\bar
v)) \, ds \right)^{1/2} && \text{if} \, v < v_\infty.
\endalignedat 
\endcases
$$
Since we need $(v_B, u_B)$ to be on this curve, we obtain that the set
of $(v_\infty, u_\infty)$ so that (4.10) has a solution lies on the
curve 
$$
u_\infty = \cases 
\alignedat2 
& u_B - \left(\int^{v_B}_{v_\infty} 2 \, (\sigma(s)-\sigma(\bar
v)) \, ds\right)^{1/2} \quad && \text{if} \ v_\infty< v_B \\
& u_B + \left(\int_{v_\infty}^{v_B} 2 \, (\sigma(s) -\sigma(\bar
v)) \, ds \right)^{1/2} && \text{if} \ v_\infty > v_B, 
\endalignedat 
\endcases
\tag 4.14
$$
where $(v_B, u_B)$ is fixed. 
\vskip.3cm 

Let us now turn to the Lax Friedrichs scheme. For the system (4.8), 
the discrete boundary layer equation is 
$$
H(v(y), v(y+1), u_\infty,  v_\infty) \equiv 
\left(
\aligned
& \lambda (u(y+1) + u(y)) + v(y+1) - v(y) -2 \, \lambda \, u_\infty) \\
& \lambda (\sigma(v(y+1)) + \sigma(v(y)) + u(y+1) - u(y) - 2
\, \lambda \, \sigma(v_\infty)
\endaligned
\right) = 0, 
\tag 4.15 
$$
$$
(v, u) (0) = (u_B, u_0), \qquad (v, u) (\infty) = (v_\infty, u_\infty).
$$
Here the eigenvalues of (4.8) are 
$$
\lambda_2(v_\infty, u_\infty) = - \lambda_1(v_\infty, u_\infty) 
                              =  \sigma'(v_\infty)^{1/2}
$$ 
and, with the notations of Section 3, 
$$
a_1(v_\infty, u_\infty) = \frac{1 -\lambda \sigma'(v_\infty)^{1/2}}
                          {1 +\lambda \sigma'(v_\infty)^{1/2}}, 
\qquad 
a_2(v_\infty, u_\infty) 
= \frac{1 + \lambda \sigma'(v_\infty)^{1/2}}
     {1 - \lambda\sigma'(v_\infty)^{1/2}}.
$$
Thus $0 < a_1 (v_\infty, u_\infty) <1$, $a_2(v_\infty, u_\infty)> 1$. 
By the analysis of Section 3, it follows that
the set of $(v_\infty, u_\infty)$ near $(v_B, u_B)$
for which (4.15) has a solution lie on a curve passing through
$(v_B, u_B)$. 
\vskip.3cm

\noindent{\bf 4.5 Eulerian Isentropic Gas Dynamics.} \, 
We now consider the isentropic approximation to the compressible
Euler system.  The system is composed of the two conservation
laws for the mass and the momentum of a gas \cite{\refCourantFriedrichs}:
$$
\aligned 
& \del_t\rho + \del_x(\rho u) = 0, \\ 
& \del_t(\rho u) + \del_x(\rho u^2 +p(\rho)) = 0. 
\endaligned
\tag 4.16 
$$
The main unknowns are the specific density $\rho$ and  
the velocity $u$. The pressure is a function of the density
and, for simplicity, we shall restrict to a 
polytropic perfect gas:
$$
p(\rho) = \rho^\gamma, \qquad \gamma \in (1, \infty). 
\tag 4.17 
$$
We consider the boundary layer equation generated by the 
vanishing viscosity method with $B(u) = I$:
$$
\aligned
& \partial_y (\rho u) = \partial^2_y \rho \\
& \partial_y (\rho u^2 + p(\rho)) = \partial^2_y u \\
& \rho(0) = \rho_B, \quad u(0) = u_B, \quad \rho(\infty) = \rho_\infty,
\quad  u(\infty) = u_\infty. 
\endaligned 
\tag 4.18
$$ 
Integrating the ODE'S and using the boundary condition at infinity, 
we get
$$
\aligned 
& \partial_y \rho = \rho u - \rho_\infty \, u_\infty \\
& \partial_y u = \rho u^2 + p(\rho) - \rho_\infty u_\infty^2 
- p(\rho_\infty) \\
& \rho(0) = \rho_B, \quad u(0) = u_B, \quad 
\rho(\infty) = \rho_\infty, u(\infty) = u_\infty. 
\endaligned
\tag 4.19
$$
The eigenvalues of the matrix obtained by linearizing the R.H.S.
of (4.19) around $(\rho_\infty, u_\infty)$ are
$$
\lambda_1(\rho_\infty,u_\infty) = u_\infty - c(\rho_\infty),
\lambda_2(\rho_\infty, u_\infty) = u_\infty + c(\rho_\infty) 
\tag 4.20
$$ 
where $c^2(\rho_\infty) = p'(\rho)$.
We have to distinguish between five different cases.  We define the
following regions in $(\rho,u)$--plane: 
$$
\aligned 
&\Omega_I = \big\{(\rho, u): u - c(\rho) < 0, u + c(\rho) < 0 \big\} \\
&\Omega_{II} = \big\{(\rho, u): u - c(\rho) < 0, u + c(\rho) =  0 \big\} \\
&\Omega_{III} = \big\{(\rho, u): u - c(\rho) < 0, u + c(\rho) > 0 \big\} \\
&\Omega_{IV}  = \big\{(\rho, u): u - c(\rho) = 0, u + c(\rho) > 0 \big\} \\
&\Omega_V= \big\{(\rho, u): u - c(\rho) > 0, u + c(\rho) > 0 \big\} 
\endaligned
\tag 4.21
$$
Thus in $\Omega_I$ both eigenvalues are negative, whereas in
$\Omega_{II}$ one has $\lambda_1 <0, \lambda_2 = 0$. In $\Omega_{III}$, one has
$\lambda_1 < 0, \lambda_2 >0$, wheras in $\Omega_{IV}$, one has 
$\lambda_1 = 0,
\lambda_2 >0$ and in $\Omega_V$, $\lambda_1 >0$ and $\lambda_2 >0$.  
Following the analysis that we did for the proof of Theorem
3.2, it is not hard to get the following local result. 
\vskip.3cm

\noindent{\it Case 1 : } $(\rho_B, u_B) \in \Omega_I$.  In this case the
set of  $(\rho_\infty, u_\infty)$ close to $(\rho_B, u_B)$ for which
(4.19) has a solution is an open neighborhood of $(\rho_B, u_B)$. 
\vskip.3cm 

\noindent{\it Case 2 : }  $(\rho_B, u_B) \in \Omega_{II}$.  In this case the
set of $(\rho_\infty, u_\infty)$ close to $(\rho_B, u_B)$ for which
(4.19) has a solution is a union of a two-dimensional region $U$ in
$\Omega_I$ and a curve in $\Omega_{III}$ through $(\rho_B, u_B)$
intersecting $U$.  
\vskip.3cm 

\noindent{\it Case 3 : } $(\rho_B, u_B) \in \Omega_{III}$.  In this case
the set of states $(\rho_\infty, u_\infty)$ close to $(\rho_B, u_B)$
for which (4.19) has a solution is a curve through $(\rho_B,
u_B)$ 
\vskip.3cm

\noindent{\it Case 4 : } $(\rho_B, u_B) \in \Omega_{IV}$.  In this case the
set of states $(\rho_\infty, u_\infty)$ near $(\rho_B, u_B)$ for
which (4.19) has a solution lies in a curve in $\Omega_{III}$
through $(\rho_B, u_B)$. This does not extend to $\Omega_{V}$.
\vskip.3cm 

\noindent{\it Case 5 : } $(\rho_B, u_B) \in \Omega_V$. There cannot be
any point $(\rho_\infty, u_\infty)$ in $\Omega_V$ for which (4.19)
has a solution.
\vskip.3cm

Next we consider the Lax-Friedrichs scheme. The discrete boundary layer 
problem to be solved is 
$$
\aligned 
& \lambda(\rho(y+1) v(y+1) + \rho(y) v(y)) - 2 \lambda \rho_\infty
u_\infty - (\rho(y+1) - \rho(y) ) = 0 \\
& \lambda(\rho(y+1) u(y+1)^2 + \rho(y) u(y)^2) - 2 \lambda \rho_\infty
u_\infty^2 - (\rho(y+1) v(y+1)  - \rho(y)v(y)) + \lambda(p(\rho(y+1))
+p(\rho(y)) - 2p(\rho_\infty))= 0
\endaligned  
\tag 4.22 
$$
$(\rho_B, u_B)$ given and $(\rho, u)(\infty) = (\rho_\infty, u_\infty)$. 

Given $(\rho_B, u_B)$ we determine $(\rho_\infty, u_\infty)$ close
to $(\rho_B, u_B)$ for which (4.22) has a solution.  Following
the analysis of the proof of Theorem 3.4, we get the eigenvalues 
of the linearized matrix at $(\rho_\infty, u_\infty)$ are
$$
a_1 = a_1 (\rho_\infty, u_\infty) = \frac{1 + \lambda 
\lambda_1(\rho_\infty, u_\infty)}{1 - \lambda  \lambda_1(\rho_\infty,
u_\infty)}, a_2 = a_2 (\rho_\infty, u_\infty) = \frac{1 + \lambda 
\lambda_2 (\rho_\infty, u_\infty)}{1 - \lambda  \lambda_2 (\rho_\infty, u_\infty)}
$$ 
where $\lambda_1$ and $\lambda_2$ are given by (4.20).  If
$(\rho_\infty, u_\infty) \in \Omega_I,  a_1 < 1, a_2 < 1$, 
if $(\rho_\infty, u_\infty) \in \Omega_{II},  a_1 < 1, a_2 = 1$, if $(\rho_\infty,
u_\infty) \in \Omega_{III},  a_1 < 1, a_2 > 1$, if $(\rho_\infty, 
u_\infty) \in \Omega_{IV},  a_1 = 1, a_2 >1$ and if $(\bar{\rho}, 
u_\infty) \in \Omega_V, a_1 > 1, a_2 > 1$.  It follows from the proof
of Theorem (3.4),  that if $(\rho_B, u_B) \in \Omega_I$, 
then the set of states $(\rho_\infty, u_\infty)$ near
$(\rho_B, u_B)$ for which (4.22) has a solution connecting
$(\rho_B, u_B)$ to $(\rho_\infty, u_\infty)$ is a neighborhood of $(\rho_B,
u_B)$.  If $(\rho_B,  u_B) \in \Omega_{II}$ this set is a union of
an  open set $U$  in $\Omega_I$ and a curve in
$\Omega_{III}$ through $(\rho_B, 
u_B)$ which interset $U$.  If $(\rho_B, u_B) \in \Omega_{III}$
this set of $(\rho_\infty, u_\infty)$ near $(\rho_B, u_B)$  
consists of a curve through $(\rho_B, u_B)$ and if $(\rho_B,
u_B) \in \Omega_{IV}$ this set consists of a curve in
$\Omega_{III}$  through   $(\rho_B, u_B)$.  If $(\rho_B, u_B) \in
\Omega_V $ no point  $(\rho_\infty, u_\infty) \in \Omega_V $ can be
connected by a solution of (4.22) from $(\rho_B, u_B)$. 
\vskip.3cm 

\noindent 
{\bf 4.6 Lagrangian Isentropic Gas Dynamics. } 
Finally, we consider the system of gas dynamics in Lagrangian coordinates
$$
\aligned
& \del_t v_t - \del_x u = 0, \\
& \del_t u + \del_x \left(\frac{1}{v}\right) = 0, 
\endaligned 
\tag 4.23
$$
where $u$ is the velocity and $v > 0$  is the
specific density.  The eigenvalues of (4.23) are 
$$
\lambda_1 = - \frac{1}{v} < 0, \qquad \lambda_2 = \frac{1}{v} > 0; 
\tag 4.24
$$ 
hence the boundary $x = 0$ is not characteristic.  

The purpose of this section is to provide an explicit formula
for the boundary layer set associated with the Lax-Friedrichs 
scheme. The boundary layer equation takes the form
$$
\aligned 
\lambda (u(y+1) + u(y))  - 2 \, \lambda \, u_\infty + v(y+1) - v(y) = 0 \\ 
\lambda \left(\frac{1}{v(y+1)} + \frac{1}{v(y)}\right) - 2 \, 
\frac{\lambda}{v_\infty} - u(y+1) + u(y) = 0 
\endaligned
\tag 4.25 
$$ 
with 
$$ 
(v(0), u(0)) = (v_B, u_B),  \qquad  (v, u)(\infty)  = (v_\infty, u_\infty). 
\tag 4.26
$$ 
We restrict attention to $v_B > \delta > 0$ for fixed $\delta$, 
and we determine the set of $(v_\infty, u_\infty)$ for which (4.25) has a
solution. We set 
$$
w(y) = \frac{v(y)}{\lambda} 
\tag 4.27
$$ 
so that (4.25) becomes 
$$
\frac{1}{w(y+1)}+ \frac{1}{w(y)} - u(y+1) + u(y) =
\frac{2}{w_\infty} 
$$ 
$$
w(y+1) - w(y) + u(y+1) + u(y) = 2 \, u_\infty.
$$ 
Adding the two equalities, we get 
$$
w(y+1) + \frac{1}{w(y+1)} + \frac{1}{w(y)} - w(y) + 2 u(y) =
\frac{2}{w_\infty} + 2 \,u_\infty.
$$ 
Setting 
$$
N(y) = - 2 \, u(y) + 2 \, u_\infty - \frac{1}{w(y)} +
\frac{2}{w_\infty} + w(y), 
$$ 
we obtain a quadratic equation for $w(y+1)$: 
$$
w^2(y+1) - N(y) \, w(y+1) + 1 = 0. 
\tag 4.28 
$$ 
Therefore 
$$
w(y+1) = \frac 1 2 \, \bigg(N(y) \pm \big(N(y)^2 -4\big)^{1/2}\bigg)
$$
from which we get an expression for $u(y+1)$ as well: 
$$
u(y+1) = \frac{\lambda}{2} \, N (y) \pm \frac{\lambda}{2} \, 
(N(y)^2-4)^{1/2}. 
\tag 4.29
$$ 
Observe that 
$N(\infty) = w_\infty + 1 / w_\infty$, where 
$w_\infty = v_\infty / \lambda$ and $N(\infty)^2 - 4 
= (w_\infty - 1 / w_\infty)^2$.  The product of the two roots of (4.28) is 
equal to one.  Stability requires
$w_\infty> 1$ so we choose the larger root in (4.29).  We have
finally from (4.29) and (4.25). 
$$
\aligned 
& v(y+1) = \frac{\lambda}{2} N(y) + \frac{\lambda}{2} (N(y)^2-4)^{1/2} \\ 
& u(y+1) = 2 u_\infty - u(y) + \frac{v(y)}{\lambda} -
\frac{N(y)}{2} - \frac{1}{2} (N(y)^2-4)^{1/2}. 
\endaligned
\tag 4.30 
$$ 

The Jacobian of the R.H.S. of (4.30) at $(v_\infty, v_\infty)$ is
easily seen to be 
$$
A(v_\infty, u_\infty) = \left(\aligned 
\frac{w_\infty^{2}+1}{w_\infty^2 -1} \; \; \; \; \frac{-2\lambda
w_\infty^2}{w_\infty^2 -1} \\
\frac{-2}{\lambda(w_\infty^2-1)} \;\;\;\; \frac{w_\infty^2 +1}{w_\infty^2 -1} 
\endaligned
\right), 
$$ 
whose eigenvalues are
$$
a_1 = \frac{w_\infty -1}{w_\infty + 1}, \qquad 
a_2 = \frac{w_\infty + 1}{w_\infty - 1}. 
$$ 
In terms of $v_\infty$, we have
$$
a_1 = \frac{1 - \frac{\lambda}{v_\infty}}{1 +
\frac{\lambda}{v_\infty}}, 
\qquad 
a_2 = \frac{1 + \frac{\lambda}{v_\infty}}{1-\frac{\lambda}{v_\infty}}.
$$
If the data for the Lax-Friedrichs scheme are chosen such that
the $v$ component is bounded away from zero, then so  is the
approximate solution.  Hence we can restrict attnetion to $v_\infty > \delta'$
for some $\delta' >0$.  For $\lambda$ small enough, we have
$$
0 < a_1 < 1 \quad \text {and} \quad a_2>1,
$$ 
and Theorem 3.10 applies.  We deduce that the set of all states $(v_\infty, u_\infty)$ 
near $(v_B, u_B)$ for which (4.25)-(4.26) has a solution
is a curve passing through $(v_B, u_B)$.

\subheading{5. Concluding Remarks}
Given a family of sets such as those introduced in this paper, we 
can formulate the boundary condition for the hyperbolic problem. 
When the solutions $u$ under consideration are functions of bounded 
variation, the traces exist in a strong sense and one can require that 
$$
u(0+,t) \in \EE(u_B(t)), \qquad t>0,
\tag 5.1 
$$
holds for all, except countably many, $t$. 
This type of regularity has been recently rigorously established 
by Amadori in her thesis \cite{\refAmadori}, using 
the front tracking scheme and 
the sets $\EE_{\text{Godunov}}$ 
($= \EE_{\text{Godunov}}^{\text{layer}}$ $= \EE_{\text{Godunov}}^{\text{entropy}}$). 
It would be interesting to extend \cite{\refAmadori} to the other sets 
we introduced here.

For general $L^\infty$ solutions constructed by the vanishing viscosity method, 
the boundary condition
$$
\supp \nu_{0,t} \subset \, \EE_{\text{viscosity}}^{\text{entropy}}
\tag 5.2 
$$
has been rigorously derived in Theorem 2.1. 
When the method of compensated compactness applies \cite{\refDiPernaone}, an
existence theorem for the boundary-value problem (1.1)--(1.3), (5.2)
follows immediatly from Theorem 2.1. 
Such a result is satisfactory only when the condition (5.2) yields,
for simple enough initial and boundary data at least, a well-posed problem. 
This is the case for the scalar equations and the linear systems, 
and, likely, for any system in the so-called Temple's class 
(having coinciding shock and rarefaction curves). 

In a recent preprint by Grenier and Gues,  
a scaling of the type $1/\sqrt{\eps}$ is used for linear systems of equations
to obtain a more precise description of the boundary features. 
However, as far as the 
formulation of a well-posed, limiting boundary-value problem for the hyperbolic 
equations is sought, 
the scaling $1/\eps$ we used in (2.5) happens to be sufficiently discriminating. 
The formulation based on the boundary layer equation may not be appropriate as it is
when the boundary is characteristic. On the other hand, the formulation based on entropy
inequality capture rigorously some of the features in the solution near the boundary, 
but is more difficult to work with analytically. 
Further study of the connection between the two sets for systems is 
in progress. 

\

\subheading{Acknowledgments} 
K.T.J. thanks Peter D. Lax for helpful suggestions on the
problem discussed in this paper while he was a graduate student
in the late 80's. P.G.L. is also very grateful to Peter D. Lax 
for fruitful discussions on this problem 
when he was a Courant instructor at NYU from 1990 to 1992.  
After the completion of this work (January 96), we heard 
an interesting lecture by Zhouping Xin (May 96) 
discussing also the formulation of boundary conditions for systems. 
Most of this work was carried out in July 1995 during a visit 
of P.G.L. at the Tata Institute of Fundamental Research, Bombay. 
P.G.L. has been partially supported by 
the Centre National de la Recherche Scientifique and by 
the National Science Foundation through  
NSF grants DMS-92-09326, DMS 94-01003, DMS 95-02766, 
and a Faculty Early Career Development (CAREER) award.

\

\subheading{References}

\item{[\refAmadori]} Amadori D., 
Initial-boundary value problems for nonlinear systems of conservation laws, 
Preprint SISSA Trieste, Italy, 1995. 
\item{[\refAmadoriColombo]} Amadori D. and Colombo M., 
Continuous dependence for $2 \times 2$ systems of conservation laws 
with boundary, 
Preprint SISSA Trieste, Italy, 1995. 
\item{[\refBall]} Ball J., 
A version of the fundamental theorem for Young measures, 
in Proceedings of Conf. on ``Partial Differential Equations and Continuum Models of
Phase Transitions'', 
Nice 1988, ed. D. Serre, Springer Verlag. 
\item{[\refBardosLerouxNedelec]} 
Bardos C.W., Leroux A.-Y., and Nedelec J.-C., 
First order quasilinear equations with boundary conditions, 
Comm. Part. Diff. Equa. 4 (1979), 1018--1034. 
\item{[\refBenabdallah]} Benabdallah A., 
Le p-syst\`eme dans un intervalle, 
C.R. Acad. Sc., Paris, t. 303, S\'erie I, 4 (1986), 
123--126. 
\item{[\refBenabdallahSerreone]} Benabdallah A. and Serre D., 
Probl\`emes aux limites pour des syst\`emes
hyperboliques nonlin\'eaires de deux \'equations \`a une dimension d'espace, 
C.R. Acad. Sc. Paris, S\'erie I, 305 (1987), 677--680.
\item{[\refBenabdallahSerretwo]} Benabdallah A. and Serre D., 
Probl\`emes aux limites pour des syst\`emes
hyperboliques nonlin\'eaires de deux \'equations \`a une dimension d'espace, 
Preprint 69, Univ. Paris-Nord, Villenateuse, France, 1988 (unpublished report).
\item{[\refBourdelDelormeMazet]} 
Bourdel F., Delorme P., and Mazet P.A., Convexity in hyperbolic problems, 
Application to a discontinuous Galerkin method, 
Proc. Inter. Conf. on Hyperbolic problems, Aachen (Germany), 
March 1988, Notes on Numer. Fluid Mech., Vol. 24, Vieweg 
(1989).
\item{[\refChenLeFloch]} Chen G.-Q. and LeFloch P.G., 
Entropy flux splittings for hyperbolic conservation laws.
General framework,
Comm. Pure Appl. Math. 48 (1995), 691--729. 
\item{[\refCourantFriedrichs]} Courant R. and Friedrichs K.O., 
Supersonic Flows and Shock Waves, 
Interscience Publishers Inc., New York (1948). 
\item{[\refDafermos]} Dafermos C.M., 
Hyperbolic systems of conservation laws, 
Proceedings ``Systems of Nonlinear Partial Differential Equations'', 
J.M. Ball editor, NATO Adv. Sci. Series C, 111, Dordrecht D. Reidel 
(1983), 25--70.
\item{[\refDiPernaone]} DiPerna R.J., 
Convergence of approximate solutions to conservation laws, 
Arch. Rational Mech. Anal. 82 (1983), 27--70.
\item{[\refDiPernatwo]} DiPerna R.J., 
Measure-valued solutions to conservations laws, 
Arch. Rational Mech. Anal. 88 (1985), 223-270.
\item{[\refDuboisLeFlochone]} Dubois F. and LeFloch P.G., 
Boundary condition for a system of conservation laws, 
C.R. Acad. Sc. Paris, t. 304, S\'erie 1 (1987), 75--78.
\item{[\refDuboisLeFlochtwo]} Dubois F. and LeFloch P.G., 
Boundary conditions for nonlinear hyperbolic 
systems of conservation laws, 
J. Diff. Equa. 71 (1988), 93--122.
\item{[\refDuboisLeFlochthree]} Dubois F. and LeFloch P.G., 
Boundary conditions for nonlinear hyperbolic systems of conservation laws, 
Proc. Inter. Conf. on Hyperbolic problems, Aachen (Germany), 
March 1988, Notes on Numer. Fluid Mech., Vol. 24, Vieweg 
(1989), 96-106.
\item{[\refDubrocaGallice]} Dubroca B. and Gallice G., 
R\'esultats d'existence et d'unicit\'e du probl\`eme mixte
pour des syst\`emes hyperboliques de lois de conservation monodimensionels, 
Comm. Part. Diff. Equa. 15 (1990), 59--80. 
\item{[\refGisclonone]} Gisclon M., 
Comparaison de deux perturbations singuli\`eres 
pour l'\'equations de Burgers avec conditions aux limites, 
C.R. Acad. Sc., Paris, S\'erie I, t. 316 (1993), 1011--1014. 
\item{[\refGisclontwo]} Gisclon M., 
Etude des conditions aux limites pour un syst\`eme strictement hyperbolique
via l'appro\-ximation parabolique,
J. Math. Pures Appl. 75 (1996), 485--508.  
\item{[\refGisclonSerre]} Gisclon M. and Serre D., 
Etude des conditions aux limites pour un syst\`eme strictement 
hyperbolique via l'approximation parabolique, 
C.R. Acad. Sc., Paris, S\'erie I, t. 319 (1994), 377--382. 
\item{[\refGlimm]} Glimm J.,
Solutions in the large for nonlinear hyperbolic systems of equations, 
Comm. Pure Appl. Math. 18 (1965), 697--715.
\item{[\refGoodman]} Goodman J.,
Initial boundary value problems for hyperbolic systems
of conservation laws, 
Ph.D. Thesis, Stanford University, California, 1982. 
\item{[\refHartman]} Hartman P., 
Ordinary Differential Equations, John Wiley and Sons Inc. (1964). 
\item{[\refJosephone]} Joseph K.T., 
Burgers equation in the quarter plane: a formula for the weak limit, 
Comm. Pure Appl. Math. 41 (1988), 133-149. 
\item{[\refJosephtwo]} Joseph K.T., Boundary layers in approximate solutions, 
Trans. Amer. Math. Soc. 314 (1989), 709--726. 
\item{[\refJosephLeFloch]} Joseph K.T. and LeFloch P.G., 
Boundary layers in weak solutions to hyperbolic conservation laws, 
Preprint $\# 1402$, 
Institute for Mathematics and its Applications, 
University of Minnesota, Minneapolis, May 1996. 
\item{[\refJosephVeerappa]} Joseph K.T. and Veerappa Gowda G.D., 
Explicit formula for the solution of convex conservation laws with
boundary condition,  
Duke Math. J. 62 (1991), 401-416.
\item{[\refKreiss]} Kreiss H.O., 
Initial-boundary value problems for hyperbolic systems, 
Comm. Pure Appl. Math. 23 (1970), 277--298. 
\item{[\refLaxone]} Lax P.D., 
Weak solutions of nonlinear hyperbolic equations and their numerical
computation,
Comm. Pure Appl. Math. 7 (1954), 159--193. 
\item{[\refLaxtwo]} Lax P.D., 
Hyperbolic systems of conservation laws II, 
Comm. Pure Appl. Math. 10 (1957), 537--566.
\item{[\refLaxthree]} Lax P.D., 
Hyperbolic systems of conservation laws and the mathematical 
theory of shock waves, 
Regional Conf. Series in Appl. Math. 11, SIAM, Philadelphia (1973).
\item{[\refLeFloch]} LeFloch P.G., 
Explicit formula for scalar conservation laws 
with boundary condition, 
Math. Meth. Appl. Sc. 10 (1988), 265--287.
\item{[\refLeFlochtwo]} LeFloch P.G., 
An introduction to nonclassical shocks of systems of conservation laws, 
Proc. International School on ``Theory and Numerics for Conservation Laws'', 
Freiburg, Germany, 20-24 Oct. 97, D. Kr\"oner, M. Ohlberger and 
C. Rohde ed., Lecture Notes in Computational Science and Engineering, 
Springer Verlag, 1998.  

\item{[\refLeFlochNedelec]} LeFloch P.G. and Nedelec J.-C., 
Explicit formula for weighted scalar nonlinear 
conservation laws, 
Trans. Amer. Math. Soc. 308 (1988), 667--683.
\item{[\refLeroux]} Leroux A.-Y., 
Approximation de quelques probl\`emes hyperboliques nonlin\'eaires, 
Th\`ese d'\'etat, University of Rennes, France (1979). 
\item{[\refLiYu]} Li T.-T. and Yu W.-C., 
Boundary Value Problem for Quasilinear Hyperbolic Systems,
Duke Univ. Math. Series (1985). 
\item{[\refLiuone]} Liu T.P., 
Large time behavior of initial and initial-boundary-value problems
of general systems of hyperbolic conservation laws, 
Comm. Math. Phys. 55 (1977), 163--177. 
\item{[\refLiutwo]} Liu T.P., 
Initial-boundary value problems for gas dynamics, 
Arch. Rational Mech. Anal. (1977), 137--168. 
\item{[\refLiuthree]} Liu T.P.,
The free piston problem for gas dynamics, 
J. Diff. Equa. 30 (1978), 175--191. 
\item{[\refLiufour]} Liu T.P., 
Admissible solutions of hyperbolic conservation laws, 
Mem. Amer. Math. Soc. 30 (1981). 
\item{[\refMajdaPego]} Majda A. and Pego R., 
Stable viscosity matrices for systems of conservation laws, 
J. Diff. Equa. 56, (1985) 229-262. 
\item{[\refSableTougeronone]} Sabl\'e-Tougeron M., 
M\'ethode de Glimm et probl\`eme mixte, 
Ann. Inst. Henri Poincar\'e 10 (1993), 423--443. 
\item{[\refSableTougerontwo]} Sabl\'e-Tougeron M., 
Les $N$-ondes de Lax pour le probl\`eme mixte, 
to appear. 
\item{[\refSerre]} Serre D., Personal communication, Paris, 1987. 
\item{[\refSmoller]} Smoller J., 
Shock Waves and Reaction Diffusion Equations, 
Springer-Verlag, New York, 1983.
\item{[\refSchwartz]} Schwartz J.T., 
Nonlinear Functional Analysis,
Courant Inst. Lect. Notes, 1964. 
\item{[\refSzepessy]} Szepessy A., 
Measure-valued solutions to scalar conservation laws with boundary 
conditions, 
Arch. Rational Mech. Anal. (1989), 181--193.
\item{[\refTemple]} Temple B., 
Systems of conservation laws with coinciding shock and rarefaction curves, 
in ``Nonlinear Partial Differential Equations'', J. Smoller ed., 
Amer. Math. Soc., Contemporary Math Series 17 (1982), 143--151. 
\item{[\refXin]} Xin Z.P., Communication to the authors, June 1994. 

\enddocument